\documentclass[11pt]{article}
\usepackage{graphicx,psfrag}
\usepackage{amssymb,amsfonts,amsthm}
\usepackage{amsmath,amsthm,amssymb}
\usepackage{amsfonts}
\usepackage{color}
\usepackage{graphics}

\newtheorem{theorem}{Theorem}[section]
\newtheorem{lemma}[theorem]{Lemma}
\newtheorem{cor}[theorem]{Corollary}
\newtheorem{proposition}[theorem]{Proposition}
\theoremstyle{definition}
\newtheorem{definition}[theorem]{Definition}

\newtheorem{remarks}[theorem]{Remarks}
\newtheorem{remark}[theorem]{Remark}
\newtheorem{ques}[theorem]{Question}
 
\def\mh{{\mathcal{H}}} 
\def\unu{{\mathbf{1}}} 
\def\bH{{\mathbb{H}}} 
\def\mh{{\mathcal{H}}} 
\def\calp{\mathcal{P}}   
\def\calc{\mathcal{C}}
\def\ggg{\mathcal{G}}   

\newcommand {\iv}{^{-1}}
\newcommand{\aaa}{{\mathcal A}}
\newcommand {\free}{\mathbb{F}} 
\newcommand{\pp}{{\mathcal P}}
\newcommand{\id}{{\mathrm{id}}}
\newcommand{\Con}{{\mathrm{Con}}}
\newcommand{\lio}[1]{\lim_\omega(#1)}
\newcommand{\co}[1]{\Con_\omega(#1)}
\newcommand {\N}{\mathbb{N}} 
\newcommand {\Z}{\mathbb{Z}}            
\newcommand {\R}{\mathbb{R}} 
\newcommand {\Q}{\mathbb{Q}} 
\newcommand {\C}{\mathbb{C}} 
\newcommand {\K}{\mathbb{K}} 
\newcommand {\hip}{\mathbb{H}} 
\newcommand {\sph}{\mathbb{S}} 
\newcommand {\F}{\mathbb{F}} 
\newcommand {\q}{\mathfrak q} 
\newcommand {\nn}{\mathcal N} 
\newcommand {\dist}{\mathrm{dist}} 
\newcommand{\cgs}{\mathrm{Cayley} (G,S)}
\newcommand{\cgas}{\mathrm{Cayley} (\Gamma ,S)}
\newcommand{\cgsh}{\mathrm{Cayley} (G,\, S\cup \mathcal{H})}
\newcommand{\dsh}{{\mathrm{dist}}_{S\cup \mathcal{H}}}
\newcommand{\ds}{{\mathrm{dist}}_S}

\newcommand {\cf}{\mathfrak c}

\newcommand {\pgot}{\mathfrak p}
\newcommand {\me}{\medskip}

\newcommand {\bi}{\bigskip}
\newcommand {\fn}{\footnote}
\newcommand {\Notat}{\noindent {\it{Notation}}:\thickspace } 
\begin{document}
\makeatletter
\title{ Quasi-isometry rigidity of groups}
\author{
Cornelia DRU\c{T}U\\ \\
Universit\'e de Lille I,\\
Cornelia.Drutu@math.univ-lille1.fr}
\date{ }
\maketitle

\tableofcontents

\newpage

These notes represent an updated version of the lectures given at
the summer school ``G\'eom\'etries \`a courbure n\'egative ou
nulle, groupes discrets et rigidit\'es'' held from the 14-th of
June till the 2-nd of July 2004 in Grenoble.

Many of the open questions formulated in the paper do not belong
to the author and have been asked by other people before. Note
that some questions are merely rhetorical and answered later in
the text; when a question is still open this is specified, with
the exception of Section 7, in which all questions are open.

\me

{\bf Acknowledgement.} I would like to thank the referee for
his/her comments, corrections and useful references that lead to
the improvement of the paper.

\section{Preliminaries on quasi-isometries}

\textit{Nota bene}: In order to ensure some coherence in the
exposition, many notions are not defined in the text, but in a
Dictionary at the end of the text.

\subsection{Basic definitions}

A \textit{quasi-isometric embedding} of a metric space $(X,
\dist_X)$ into a metric space $(Y,\dist_Y)$ is a map $\q : X \to
Y$ such that for every $x_1,x_2\in X$,
\begin{equation}\label{qi}
\frac{1}{L}\dist_X (x_1,x_2)-C \leq \dist_Y (\q(x_1), \q(x_2))\leq
L\dist_X (x_1,x_2)+C\, ,
\end{equation} for some constants $L\geq 1$ and $C\geq 0$.

If $X$ is a finite interval $[a,b]$ then $\q$ is called
\textit{quasi-geodesic (segment)}. If $a=-\infty$ or $b=+\infty$
then $\q$ is called \textit{quasi-geodesic ray}. If both
$a=-\infty$ and $b=+\infty$ then $\q$ is called
\textit{quasi-geodesic line}. The same names are used for the
image of $\q$.

\me

If moreover $Y$ is contained in the $C$--tubular neighborhood of
$\q (X)$ then $\q$ is called \textit{a quasi-isometry}. In this
case there exists $\bar{\q} :Y \to X$ quasi-isometry such that
$\bar{\q} \circ \q$ and $\q \circ \bar{\q}$ are at uniformly
bounded distance from the respective identity maps (see \cite{GH}
for a proof). We call $\bar{\q}$ \textit{quasi-converse of} $\q$.

\begin{figure}[!ht]
\centering
\includegraphics{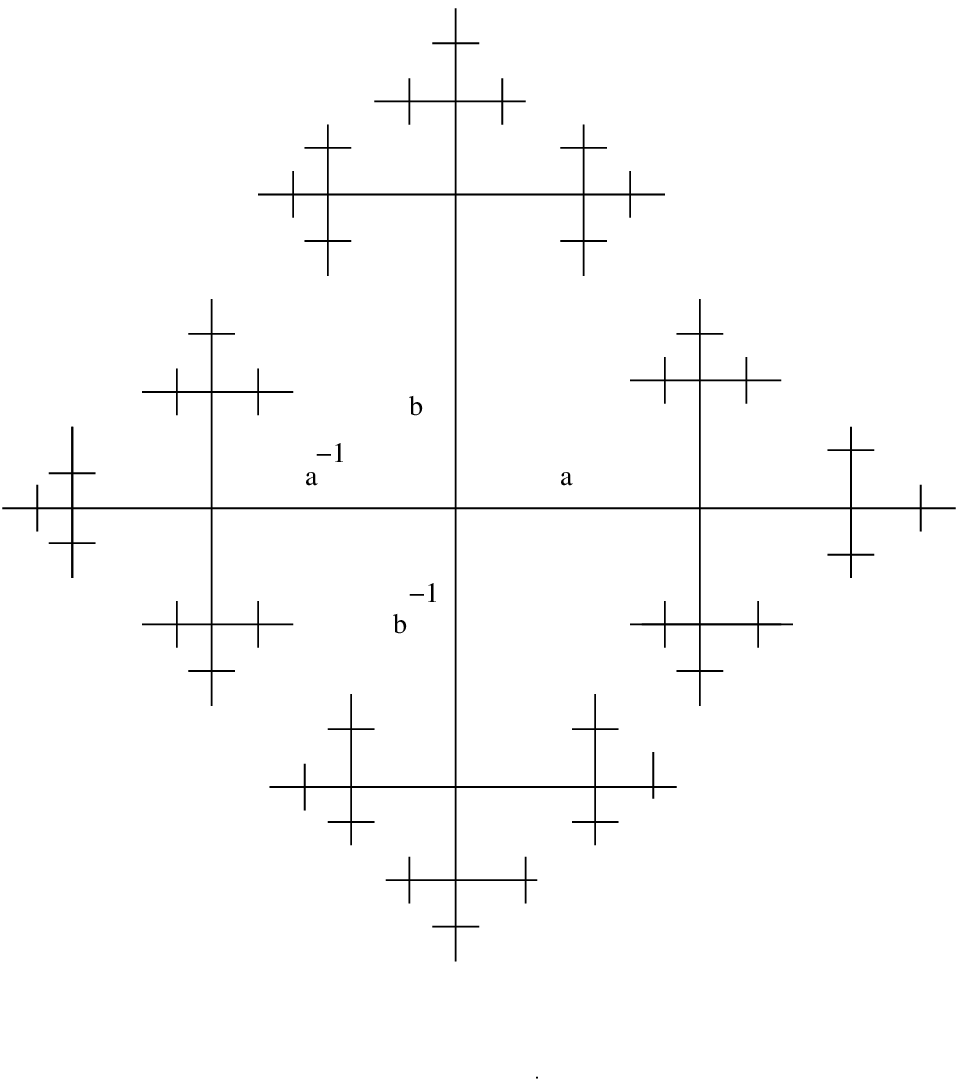}
\caption{Cayley graph of $\free_2$.} \label{fig1}
\end{figure}

The objects of study are the finitely generated groups. We first
recall how to make them into geometric objects. Given a group $G$
with a finite set of generators $S$ containing together with every
element its inverse, one can construct the \textit{Cayley graph}
$\cgs$ as follows:
\begin{itemize}
    \item its set of vertices is
$G$;
    \item every pair of elements $g_1,g_2\in G$ such that $g_1=g_2s$,
with $s\in S$, is joined by an edge. The oriented edge $(g_1,g_2)$
is labelled by~$s$.
\end{itemize}

We suppose that every edge has length $1$ and we endow $\cgs$ with
the length metric. Its restriction to $G$ is called \textit{the
word metric associated to $S$} and it is denoted by $\dist_S$. See
Figure \ref{fig1} for the Cayley graph of the free group of rank
two $\free_2=\langle a,b\rangle$.

\begin{remark}\label{cayinf}
A Cayley graph can be constructed also for an infinite set of
generators. In this case the graph has infinite valence in each
point.
\end{remark}

\begin{figure}[!ht]
\centering
\includegraphics{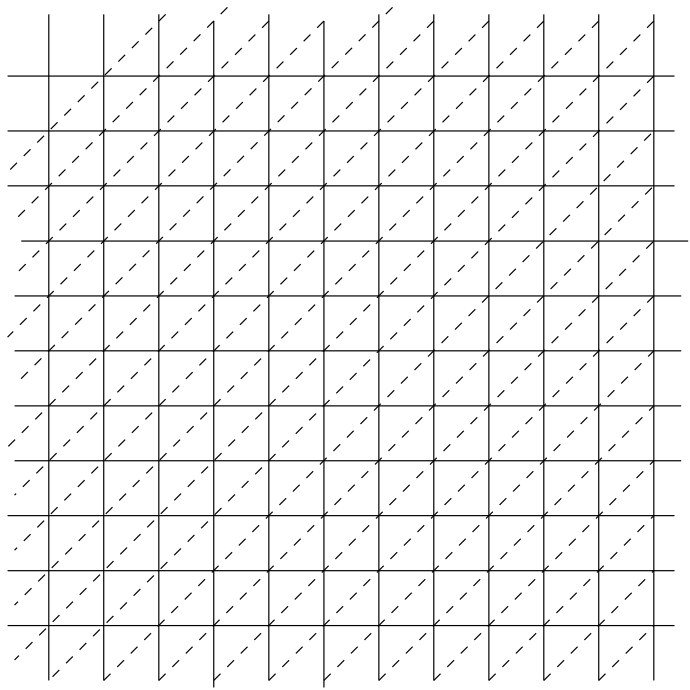}
\caption{Cayley graph of $\Z^2$.} \label{fig2}
\end{figure}

We note that if $S$ and $\bar{S}$ are two finite generating sets
of $G$ then $\dist_S$ and $\dist_{\bar{S}}$ are bi-Lipschitz
equivalent. In Figure \ref{fig2} are represented the Cayley graph
of $\Z^2$ with set of generators $\{ (1,0), (0,1)\}$ and the
Cayley graph of $\Z^2$ with set of generators $\{ (1,0), (1,1)\}$.

\me

\subsection{Examples of quasi-isometries}\label{exqi}

\begin{enumerate}
    \item The main example, which partly justifies
    the interest in quasi-isometries, is the following. Given $M$ a
    compact Riemannian manifold, let $\widetilde{M}$ be its
    universal covering and let $\pi_1(M)$ be its fundamental group. The group $\pi_1(M)$ is finitely
    generated, in fact even finitely presented \cite[Corollary I.8.11, p.137]{BH}.

     The metric space $\widetilde{M}$ with the Riemannian metric is quasi-isometric to
    $\pi_1(M)$ with some word metric. This can be clearly seen in
    the case when $M$ is the $n$-dimensional flat torus
    $\mathbb{T}^n$. In this case $\widetilde{M}$ is $\R^n$ and
    $\pi_1(M)$ is $\Z^n$. They are quasi-isometric, as $\R^n$ is a
    thickening of $\Z^n$.
    \item More generally, if a group $\Gamma $ acts properly discontinuously
    and with compact quotient by isometries on a complete locally compact length metric space $(X, \dist_\ell)$
    then $\Gamma $ is finitely generated \cite[Theorem I.8.10, p. 135]{BH} and
    $\Gamma$ endowed with any word metric is quasi-isometric to
    $(X\, ,\, \dist_\ell)$.

    Consequently two groups acting as above on the same length metric space are
    quasi-isometric.

    \item Given a finitely generated group $G$ and a finite
    index subgroup $G_1$ in it, $G$ and $G_1$ endowed with
    arbitrary word metrics are quasi-isometric. This may be seen
    as a particular case of the previous example, with $\Gamma =
    G_1$ and $X$ a Cayley graph of $G$.

    In terms of Riemannian manifolds, if $M$ is a finite covering of
    $ N$ then $\pi_1(M)$ and $\pi_1(N)$ are quasi-isometric.
    \item Given a finite
    normal subgroup $N$ in a finitely generated group $G$, $G$ and $G/N$ (both endowed with
    arbitrary word metrics) are quasi-isometric.

    Thus, in arguments where we
    study behavior of groups with respect to quasi-isometry, we can
    always replace a group with a finite index subgroup or with a
    quotient by a finite normal subgroups.
    \item All non-Abelian free groups of finite rank are quasi-isometric to
    each other. This follows from the fact that the Cayley graph
    of the free group of rank $n$ with respect to a set of $n$
    generators and their inverses is the regular simplicial tree of valence $2n$.

    Now all regular
    simplicial trees of valence at least $3$
    are
    quasi-isometric. We denote by $\mathcal{T}_{k}$ the
    regular simplicial tree of valence $k$ and
    we show that $\mathcal{T}_3$ is quasi-isometric to $\mathcal{T}_k$ for every $k\geq 4$.

    We define the map $\q :\mathcal{T}_3 \to \mathcal{T}_k$ as
      in Figure \ref{fig3}, sending all edges drawn in thin lines
      isometrically onto edges and all paths of length $k-3$ drawn in thick lines
      onto one vertex. The map $\q$ thus defined is surjective and
      it satisfies the inequality
      $$
\frac{1}{k-2}\, \dist (x,y) - 1 \leq \dist (\q(x),\q(y))\leq \dist
(x,y) \, .
      $$

\begin{figure}[!ht]
\centering
\includegraphics{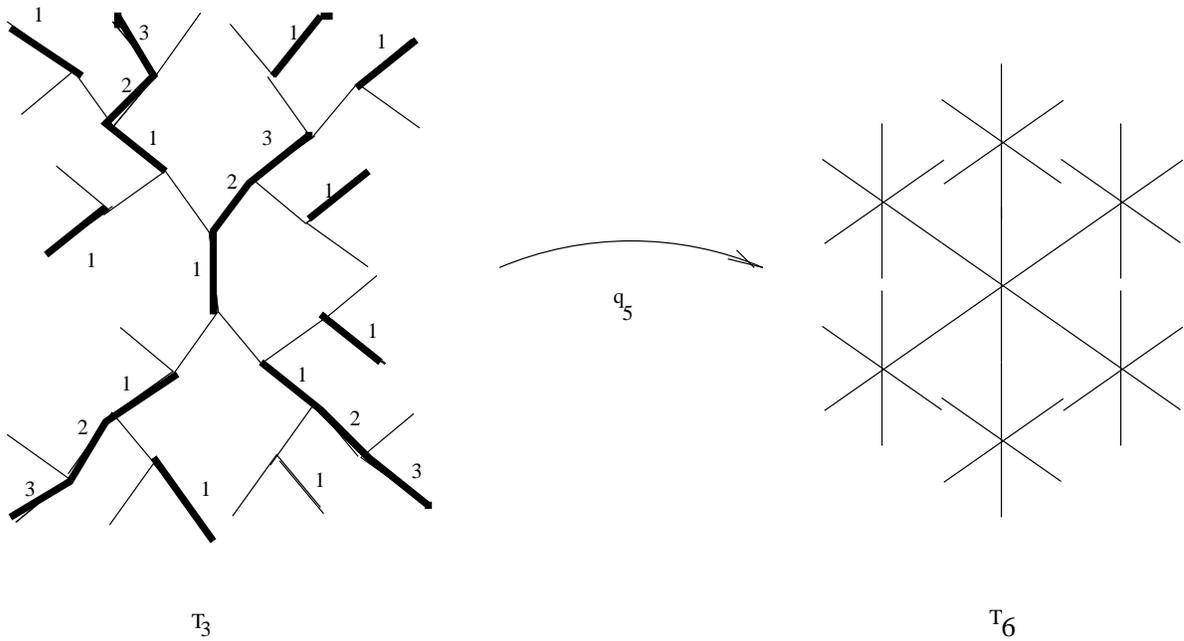}
\caption{All regular simplicial trees are quasi-isometric.}
\label{fig3}
\end{figure}

      \item Let $M$ be a non-compact hyperbolic two-dimensional orbifold of finite area. This is the same thing
      as saying that $M=\Gamma \backslash \hip^2_\R $, where $\Gamma$ is a
      discrete subgroup of $PSL_2(\R )$ with fundamental domain of
      finite area.

\noindent{\it{Nota bene}}: We assume that all the actions of
groups by isometries on spaces are to the left, as in the
particular case when the space is the Cayley graph. This means
that the quotient will be always taken to the left. We feel sorry
for all people which are accustomed to the quotients to the right.

We can apply the following result.

\begin{lemma}[Selberg's Lemma]
A finitely generated group which is linear over a field of
characteristic zero has a torsion free subgroup of finite index.
\end{lemma}

We recall that \textit{torsion free group} means a group which
does not have finite non-trivial subgroups. For an elementary
proof of Selberg's Lemma see \cite{Al}.

We conclude that $\Gamma$ has a finite index subgroup $\Gamma_1$
which is torsion free. It follows that $N = \Gamma_1 \backslash
\hip^2_\R $ is a hyperbolic surface which is a finite covering of
$M$, hence it is of finite area but non-compact. On the other
hand, it is known that the fundamental group of such a surface is
a free group of finite rank (see for instance \cite{Ma}).

\textit{Conclusion}: The fundamental groups of all hyperbolic
two-dimensional orbifolds are quasi-isometric to each other.
\end{enumerate}

\medskip

At this point one may start thinking that the quasi-isometry is
too weak a relation for groups, and that it does not distinguish
too well between groups with different algebraic structure. It
goes without saying that we are discussing here only infinite
finitely generated groups, because we need a word metric and
because finite groups are all quasi-isometric to the trivial
group.

We can begin by asking if the result in Example 6 is true for any
rank one symmetric space.

\me

\begin{ques}\label{qs}
Given $M$ and $N$ orbifolds of finite volume covered by the same
rank one symmetric space, is it true that $\pi_1(M)$ and
$\pi_1(N)$ are quasi-isometric ?
\end{ques}

\me

It is true if $N$ is obtained from $M$ by means of a sequence of
operations obviously leaving the fundamental group $\Gamma
=\pi_1(M)$ in the same quasi-isometry class :
\begin{itemize}
    \item going up or down a finite covering, which at the level
    of fundamental groups means changing $\Gamma $ with a finite index
    subgroup or a finite extension;
    \item replacing a manifold with another one isometric to it,
    which at the level of groups means changing $\Gamma $ with a
    conjugate $\Gamma^g$, where $g$ is an isometry of the
    universal covering.
\end{itemize}

Above we have used the following

\me

\noindent \textit{Notation}: For $A$ an element or a subgroup in a
group $G$ and $g\in G$, we denote by $A^g$ its image $gA g^{-1}$
under conjugacy by $g$.

\me

In the sequel we also use the following

\me

\noindent \textit{Convention}: In a group $G$ we denote its
neutral element by $\id$ if we consider an action of the group on
a space, and by $1$ otherwise.

\section{Rigidity of non-uniform rank one lattices}

It turns out that the answer to the Question \ref{qs} is ``very
much negative'', so to speak, that is: apart from the exceptional
case of two dimensional hyperbolic orbifolds, in the other cases
finite volume rank one locally symmetric spaces which are not
compact have quasi-isometric fundamental groups if and only if the
locally symmetric spaces are obtained one from the other by means
of a sequence of three of the operations described previously.
More precisely, the theorem below, formulated in terms of groups,
holds.

\subsection{Theorems of Richard Schwartz}

We recall that a discrete group of isometries $\Gamma$ of a
symmetric space $X$ such that $\Gamma \backslash X$ has finite
volume is called a \textit{lattice}. If $\Gamma \backslash X$ is
compact, the lattice is called \textit{uniform}, otherwise it is
called \textit{non-uniform}.

\begin{theorem}[R. Schwartz, \cite{Sch}]\label{RSch}
\begin{itemize}
    \item[(1)](\textit{quasi-isometric lattices}) Let $G_i$ be a non-uniform lattice of isometries of
    the rank one symmetric space $X_i,\, i=1,2$. Suppose that $G_1$
    is quasi-isometric to $G_2$. Then $X_1=X_2=X$ and one of the
    following holds:
    \begin{itemize}
        \item[(a)] $X=\hip^2_\R$;
        \item[(b)] there exists an isometry $g$ of $X$ such that
        $G_1^{g}\cap G_2$ has finite index both in $G_1^{g}$ and in $G_2$.
    \end{itemize}
    \item[(2)](\textit{finitely generated groups quasi-isometric to
    lattices}) Let $\Lambda$ be a finitely generated group and let $G$ be
     a non-uniform lattice of isometries of a rank one symmetric space $X\neq
     \hip^2_\R$. If $\Lambda $ is quasi-isometric to $G$ then there
     exists a non-uniform lattice $G_1$ of isometries of $X$ and a finite group $F$ such that one has the following exact
     sequence:
     $$
0\to F\to \Lambda \to G_1 \to 0 \, .
     $$
\end{itemize}
\end{theorem}

The notion of commensurability is recalled in Section \ref{dic}.
The particular case of commensurability described in Theorem
\ref{RSch}, (1), (b), means that the locally symmetric spaces
$G_1\backslash X$ and $G_2\backslash X$ have isometric finite
coverings.

We note that Theorem \ref{RSch} is in some sense a much stronger
result than Mostow Theorem. Mostow Theorem requires the
isomorphism of fundamental groups - which is an algebraic relation
between groups, also implying their quasi-isometry. Theorem
\ref{RSch} only requires that the groups are quasi-isometric,
which is a relation between ``large scale geometries'' of the two
groups, and has \textit{a priori} nothing to do with the algebraic
structure of the groups.

Since Mostow rigidity holds for all kinds of lattices, a first
natural question to ask is:

\medskip

\begin{ques}\label{schu}
 Are the two statements in Theorem \ref{RSch}
also true for uniform lattices~?
\end{ques}

\medskip

Concerning statement (1), the following can be said: given $G_i$
uniform lattices of isometries of
    the rank one symmetric spaces $X_i,\, i=1,2$, $G_1$ quasi-isometric to $G_2$ implies that
    $X_1=X_2=X$.

   Now one can ask if in case $X\neq \hip_\R^2$ there exists an isometry $g$ of $X$ such that
        $G_1^{g}\cap G_2$ has finite index both in $G_1^{g}$ and in $G_2$ ? In other words is it true that
        \textit{all} uniform lattices of isometries of the same
        rank one symmetric space $X\neq \hip_\R^2$ are
        commensurable ?

        A weaker variant of the previous question is whether all
        \textit{arithmetic} uniform lattices of isometries of  $X\neq \hip_\R^2$ are commensurable.

        The answer to both questions is negative, as shown by the
        following counter-example.

        \me

\textit{Counter-example}:

\me

All the details for the statements below can be found in
\cite{GPS}.

Let $Q$ be a quadratic form of the type $\sqrt{2} x_{n+1}^2 -a_1
x_1^2 -\cdots - a_n x_n^2$, where $a_i$ are positive rational
numbers. The set
$$
\hip_Q =\{ (x_1,\dots , x_{n+1})\in \R^{n+1}\mid Q(x_1,\dots ,
x_{n+1})=1\, ,\, x_{n+1}>0\}
$$ is a model of the hyperbolic $n$-dimensional space. Its group
of isometries is $SO_{\mathrm{Id}}(Q)$, the connected component
containing the identity of the stabilizer of the form $Q$ in
$SL(n+1,\R)$. The discrete subgroup $G_Q = SO_{\mathrm{Id}}(Q)
\cap SL(n+1, \Z (\sqrt{2}))$ is a uniform lattice. Now if two such
lattices $G_{Q_1}$ and $G_{Q_2}$ are commensurable then there
exist $g\in GL(n+1,\Q [\sqrt{2}])$ and $\lambda \in \Q [\sqrt{2}]$
such that $Q_1 \circ g = \lambda Q_2$. In particular, if $n$ is
odd then the ratio between the discriminant of $Q_1$ and the
discriminant of $Q_2$ is a square in $\Q [\sqrt{2}]$. It now
suffices to take two forms such that this is not possible, for
instance (like in \cite{GPS}):
$$
Q_1=\sqrt{2} x_{n+1}^2 -x_1^2 - x_2^2 -\cdots - x_n^2 \mbox{ and
}Q_2=\sqrt{2} x_{n+1}^2 -3x_1^2 - x_2^2 -\cdots - x_n^2\, .
$$

\me

Statement (2) of Theorem \ref{RSch}, on the other hand, also holds
for uniform lattices. See the discussion in the beginning of
Section \ref{cgc}.

\me

 A main step in the proof of Theorem \ref{RSch} is the following rigidity
result, interesting
 by itself.

 \begin{theorem}[Rigidity Theorem \cite{Sch}]\label{rig}
Let $\Gamma $ and $G$ be two non-uniform lattices of isometries of
$\hip^n_\F \neq \hip^2_\R \, ,\, \F= \R $ or $\C$. An
$(L,C)$--quasi-isometry $\q$ between $\Gamma $ and $G$ is at
finite distance from an isometry $g$ of $\hip^n_\F$ with the
property that $\Gamma^g\cap
    G$ has finite index in both $\Gamma^g$ and $G$.
 \end{theorem}

 The meaning of the statement ``\textit{$\q$ is at finite distance from $g$}'' is the following:

\begin{center}
For every compact $\mathcal{K}$ in $G \backslash \hip^n_\F$ there
exists $D=D(L,C,\mathcal{K},\Gamma ,G)$ such that for every $x_0$
with $G x_0 \in \mathcal{K}$, one has
$$
\dist (\q (\gamma )x_0\, ,\, g\gamma x_0)\leq D,\: \forall \gamma
\in \Gamma \, .
$$
\end{center}

As it is, it does not look very enlightening. We shall come back
to this statement in Section \ref{pfrig}, after recalling what is
the structure of finite volume real hyperbolic manifolds in
Section \ref{3man}. Also in Section \ref{pfrig} we shall give an
outline of the proof of Theorem \ref{rig} in the particular case
when $\hip_\F^n = \hip^3_\R$. All the ideas of the general proof
are already present in this particular case, and we avoid some
technical difficulties that are irrelevant in a first approach.

According to Selberg's Lemma, we may suppose, without loss of
generality, that both $\Gamma $ and $G$ are without torsion.

\subsection{Finite volume real hyperbolic manifolds}\label{3man}

Let $M$ be a finite volume real hyperbolic manifold, that is a
manifold with universal covering $\hip_\R^n$, for some $n\geq 2$.
Let $\Gamma $ be its fundamental group.

Given a point $x\in M$ denote by $r(x)$ the injectivity radius of
$M$ at $x$.

For every $\varepsilon >0$ the manifold can be decomposed into two
parts:
\begin{itemize}
    \item the $\varepsilon $-\textit{thick part}: $M_{\geq \varepsilon }=\{ x\in M \mid r(x) \geq \varepsilon
    \}$;
    \item  the $\varepsilon $-\textit{thin part}: $M_{< \varepsilon }=\{ x\in M \mid r(x) < \varepsilon
    \}$.
\end{itemize}

The following theorem describes the structure of $M$. We refer to
\cite[$\S 4.5$]{Th} for details.

\begin{theorem}\label{th}
\begin{itemize}
    \item[(1)] There exists a universal constant $\varepsilon_0
    =\varepsilon_0
(n)>0$ such that for every complete manifold $M$ of universal
covering $\hip_\R^n$ and of fundamental group $\Gamma $, the
$\varepsilon_0 $-thin part $M_{< \varepsilon_0}$ is a disjoint
union of
\begin{itemize}
    \item tubular neighborhoods of short closed geodesics;
    \item neighborhoods of cusps, that is sets of the form $\Gamma_\alpha \backslash Hbo(\alpha
    ) $, where $Hbo (\alpha )$ is an open horoball of basepoint
    $\alpha \in \partial_\infty \hip_\R^n$ and $\Gamma_\alpha $ is the stabilizer of $\alpha $ in $\Gamma
    $.
\end{itemize}
    \item[(2)] A complete hyperbolic manifold $M$ has finite volume if and only if for every $\varepsilon >0$
     the $\varepsilon $-thick part $M_{\geq \varepsilon }$ is
     compact.

    Note that the fact that $M_{\geq \varepsilon_0 }$ is
     compact implies that
    \begin{itemize}
        \item  $M_{< \varepsilon_0}$ has finitely many components;
        \item  for every neighborhood of a cusp,  $\Gamma_\alpha \backslash Hbo(\alpha
    ) $, its boundary $\Gamma_\alpha \backslash H(\alpha )$, where $H(\alpha)$ is
    the boundary horosphere of $Hbo(\alpha )$, is compact.
    \end{itemize}
\end{itemize}

\end{theorem}

Let now $M$ be a finite volume real hyperbolic manifold and
$\Gamma =\pi_1(M)$. Consider the finite set of tubular
neighborhoods of cusps
$$
\{ \Gamma_{\alpha_i}\backslash Hbo(\alpha_i
    ) \mid i\in \{ 1,2,\dots ,m \} \}\, .
$$

According to Theorem \ref{th}, the set
$$
M_0=M\setminus \coprod_{i=1}^m \Gamma_{\alpha_i}\backslash
 Hbo(\alpha_i) $$ is compact. The pre-image of each cusp $\Gamma_{\alpha_i}\backslash Hbo(\alpha_i
    )$ is the $\Gamma$--orbit of $Hbo(\alpha_i
    )$ and the open horoballs composing this orbit are pairwise
    disjoint. Thus, the space
    $$
X_0=\hip^3_\R \setminus \coprod_{i=1}^m \coprod_{\gamma \in \Gamma
/\Gamma_{\alpha_i}} \gamma Hbo(\alpha_i
    )
    $$
     satisfies $\Gamma \backslash X_0 =M_0$.

\begin{figure}[!ht]
\centering
\includegraphics{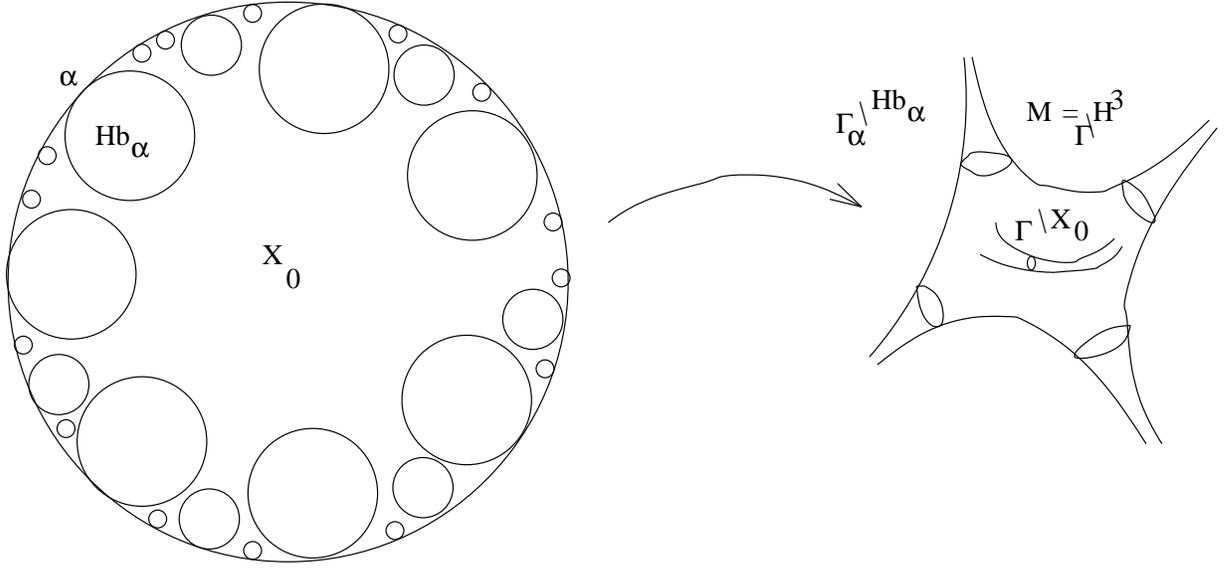}
\caption{Finite volume hyperbolic manifold, space $X_0$.}
\label{fig4}
\end{figure}

\begin{remarks}
The group $\Gamma$ endowed with a word metric $\dist_w$ is
quasi-isometric to the space $X_0$ with the length metric
$\dist_\ell$ according to the example (2) of quasi-isometry given
in Section \ref{exqi}. In particular, since for every $x_0\in
X_0$, $\Gamma x_0$ is a net in
    $X_0$, we also have that $(\Gamma \, ,\, \dist_w)$ is
quasi-isometric to $(\Gamma x_0\, ,\, \dist_\ell)$.
\end{remarks}

\subsection{Proof of Theorem \ref{rig}}\label{pfrig}

According to Section \ref{3man}, there exists $X_0$ complementary
set in $\hip^3_\R$ of countably many pairwise disjoint open
horoballs such that $\Gamma \backslash X_0$ is compact.
Consequently, $\Gamma$ with any word metric is quasi-isometric to
$X_0$ with the length metric. Similarly, one can associate to $G$
a complementary set $Y_0$ of countably many pairwise disjoint open
horoballs such that $G\backslash Y_0 $ is compact and such that
$G$ with a word metric is quasi-isometric to $Y_0$ with the length
metric.

The quasi-isometry $\q :\Gamma \to G$ induces a quasi-isometry
between $(\Gamma x_0 , \dist_\ell)$ and $(G y_0,\dist_\ell )$, for
every $x_0\in X_0\, ,\, y_0\in Y_0$, hence also between $X_0$ and
$Y_0$ (each quasi-isometry having different parameters $L$ and
$C$). For simplicity, we denote all these quasi-isometries by $\q$
and all their constants by $L$ and $C$.

In these terms, the conclusion of Theorem \ref{rig} means that
$\q$ seen as a quasi-isometry between $X_0$ and $Y_0$ is at
distance at most $D$ from the restriction to $X_0$ of an isometry
$g$ in $Comm (\Gamma ,G)$, where $D=D(L,C,X_0,Y_0)$.

We now give an outline of proof of Theorem \ref{rig}.

\me

\textsc{Step 1}. The following general statement holds.

\begin{lemma}[Quasi-Flat Lemma \cite{Sch}, $\S 3.2$]\label{qiflats}
Let $\Gamma$ be a non-uniform lattice of isometries of
$\hip^3_\R$. For every $L\geq 1$ and $C\geq 0$ there exists
$M=M(L,C)$ such that every quasi-isometric embedding
$$
\q :\Z^2 \to \Gamma
$$ has its range in $\nn_M (\gamma \Gamma_\alpha )$, where
$\Gamma_\alpha$ is a cusp group and $\gamma \in \Gamma$.
\end{lemma}

Let us apply Lemma \ref{qiflats} to the $(L,C)$--quasi-isometry
$\q$ from $\Gamma $ to $G$, and to its quasi-converse $\bar{\q} :
G \to \Gamma $.
\begin{itemize}
    \item For every $\gamma \in \Gamma $ and $\alpha \in
    \partial_\infty \hip^3_\R
    $ corresponding to a cusp of $\Gamma \backslash \hip^3_\R  $,
    $$
    \q (\gamma \Gamma_\alpha )\subset \nn_M (g G_\beta )\, ,
    $$ for some $g\in G$ and $\beta \in
    \partial_\infty \hip^3_\R
    $ corresponding to a cusp of $G\backslash \hip^3_\R $;
    \item For every $g\in G$ and $\beta \in
    \partial_\infty \hip^3_\R
    $ corresponding to a cusp of $G \backslash \hip^3_\R $,
     $$
    \bar{\q} (g G_\beta )\subset \nn_M (\gamma' \Gamma_{\alpha'} )\, ,
    $$ for some $\gamma' \in \Gamma$ and $\alpha' \in
    \partial_\infty \hip^3_\R
    $ corresponding to a cusp of $\Gamma \backslash \hip^3_\R
    $.
\end{itemize}

Combining both and noticing that if the left coset $\gamma
\Gamma_\alpha$ is contained in the tubular neighborhood of another
left coset $\gamma' \Gamma_{\alpha'}$ then the two coincide, a
bijection is obtained between left cosets $\gamma \Gamma_\alpha$
and left cosets $g G_\beta$ such that
\begin{equation}\label{corr}
    \dist_H (\q (\gamma \Gamma_\alpha ) , g G_\beta )\leq M'\, .
\end{equation}

Here $\dist_H $ denotes the Hausdorff distance (see the Dictionary
for a definition).

\begin{figure}[!ht]
\centering
\includegraphics{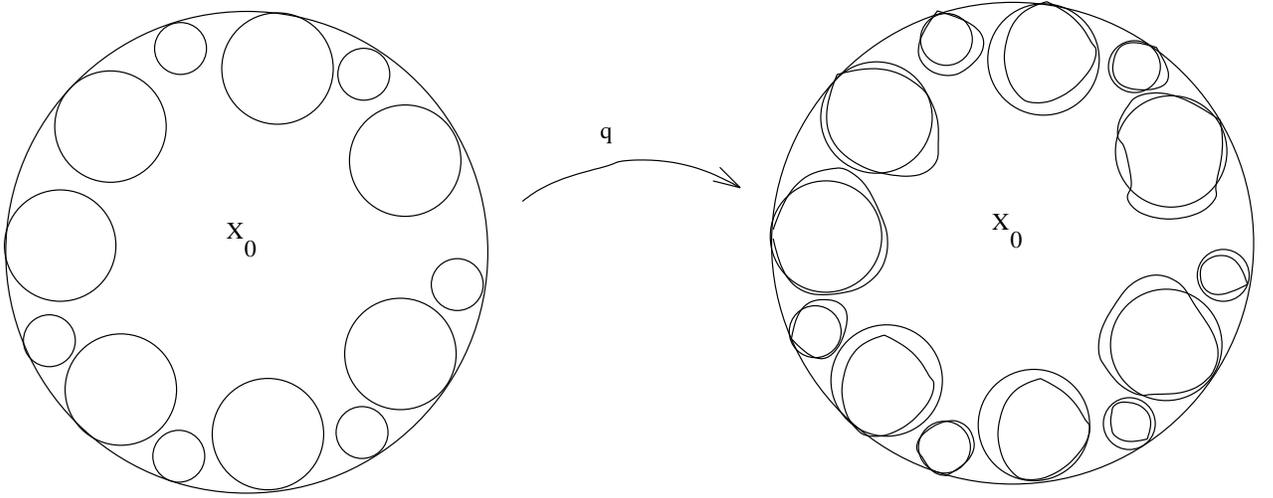}
\caption{Quasi-isometric embeddings of $\Z^2$.} \label{fig5}
\end{figure}

The situation is represented in Figure \ref{fig5}.

\me

\textsc{Step 2}. The map $\q $ seen as a quasi-isometry between
$\Gamma x_0$, net in $X_0$ and $Gy_0$, net in $Y_0$, is extended
to a quasi-isometry $\q_e$ between a net $N_1$ in $\hip_\R^3$ and
a net $N_2$ in $\hip_\R^3$. This is done horoball by horoball. Let
$\gamma \Gamma_\alpha$ and $g G_\beta$ satisfying (\ref{corr}).
Let $\gamma Hb_\alpha$ and $g Hb_\beta$ be the corresponding
horoballs. We divide each of them into strips of constant width,
by means of countably many horospheres. We note that $\gamma
\Gamma_\alpha x_0$ is $\delta$-separated, and that the horosphere
$\gamma H_\alpha $ is contained in $\nn_\epsilon (\gamma
\Gamma_\alpha x_0 )$, for some $\delta >0$ and $\epsilon >0$.

We project $\gamma \Gamma_\alpha x_0$ onto the first horosphere
$H_1$. We get a $(\delta' \, , \, \epsilon')$--net,
$\mathcal{N}_1^{(1)}$, for some $\delta'<\delta$ and $\epsilon'
<\epsilon$. We choose a maximal $\delta$-separated subset
$N_1^{(1)}$ in $\mathcal{N}_1^{(1)}$, hence a $(\delta,
\delta)$--net in $\mathcal{N}_1^{(1)}$ and a $(\delta,
\delta+\epsilon' )$--net in $H_1$. We extend $\q$ to $N_1^{(1)}$
by $\q (n_1)=\pi \circ \q \circ \pi^{-1}(n_1)$, where $\pi $
denotes the projection in each of the spaces onto the first
horosphere $H_1$.

We repeat the argument and extend $\q$ to a net in $H_2,\, H_3$,
etc. A global quasi-isometry is thus obtained. Indeed, given two
points $A\in H_n$ and $B\in H_m$, with $m\geq n$, if $B'$ is the
projection of $B$ onto $H_n$, the distance $\dist (A,B)$ is
bi-Lipschitz equivalent to $\dist (A,B')+\dist (B',B)$.

We finally obtain a quasi-isometry $\q_e$ between nets of
$\hip_\R^3$, hence a quasi-isometry of $\hip_\R^3$.

\me

\textit{Nota bene}: In the case of the Mostow rigidity theorem
also a quasi-isometry of the whole space is obtained, but it has
the extra property that it is \textit{equivariant with respect to
a given isomorphism between the two groups $\Gamma $ and $G$}.
Here, the property of equivariance is replaced by the extra
geometric information that $\q$ sends the space $X_0$ at uniformly
bounded distance from $Y_0$, by sending each boundary horosphere
at bounded distance from a boundary horosphere.

\me

The quasi-isometry $\q_e$ extends, according to the Theorem of
Efremovitch-Tikhomirova to a map between boundaries
$$
\partial \q_e : \sph^2_\infty \to \sph^2_\infty\, ,
$$ which is a quasi-conformal homeomorphism.

Next, two theorems are used.

\begin{theorem}[Rademacher-Stepanov, see \cite{LV}]\label{RS}
Every quasi-conformal homeomorphism between open sets in $\sph^2$
is differentiable almost everywhere.
\end{theorem}

\begin{theorem}[\cite{LV}]\label{mob}
A quasi-conformal homeomorphism $h:\sph^2 \to \sph^2$ whose
differential is almost everywhere a similarity is a M\"obius
transformation.
\end{theorem}

\textsc{Step 3.} It is shown that in almost every point in the set
of differentiability of $\partial \q $ the differential is a
similarity.

Let $\Omega_1$ be the set of points $\xi $ in $\sph^2_\infty$ such
that the geodesic ray $[x_0, \xi )$ returns in $X_0$ infinitely
often. The set of such $\xi $ has full Lebesgue measure in
$\sph^2_\infty$. This can be seen for instance by projecting onto
$\Gamma \backslash \hip^3_\R $ and noting that almost every
locally geodesic ray in it is equidistributed, hence it returns
infinitely often in the compact $\Gamma\backslash X_0$.

Likewise, let $\Omega_2$ be the set of points $\zeta $ in
$\sph^2_\infty$ such that the geodesic ray $[y_0, \zeta )$ returns
in $Y_0$ infinitely often.

Let $\Omega $ be the set of points in $\Omega_1 \cap \partial
\q_e^{-1} (\Omega_2 )$ in which $\partial \q_e$ is differentiable.
Let us show that in every $\xi \in \Omega$ the differential of
$\partial \q_e$ is a similarity. Denote by $\zeta$ the image
$\partial \q_e (\xi )$. Also, denote by $\xi'$ the image of $\xi$
under the geodesic symmetry of center $x_0$ and by $\zeta'$ the
image of $\zeta$ under the geodesic symmetry of center $y_0$. In
the sequel consider the two stereographic projections of
$\hip^3_\R$ sending $(\xi, \xi', x_0)$ and respectively $(\zeta,
\zeta' , y_0)$ onto $(O, \infty , (0,0,1))$. We work in the
corresponding half-space models of $\hip^3_\R$ for the domain and
the range of $\q_e$, respectively. In these models, $\partial \q_e
(O)=O$, $\partial \q_e (\infty )=\infty $ and the differential
$d_\xi
\partial \q_e$ becomes $d_O \partial \q_e$. The goal is to show that the latter is
a similarity.

Let $(x_n)$ be a sequence of points on $[x_0, \xi ) \cap X_0$
diverging to $\xi$. Let $t_n$ be the hyperbolic isometry of axis
containing $[x_0, \xi )$ such that $t_n(x_0)=x_n$.

Similarly, let $(y_n)$ be a sequence of points on $[y_0, \zeta)
\cap Y_0$ diverging to $\zeta$. Let $\tau_n$ be the hyperbolic
isometry of axis containing $[y_0, \zeta )$ such that
$\tau_n(y_0)=y_n$.

Consider the sequence of $(L,C)$--quasi-isometries
$\q_n=\tau_n^{-1} \circ \q \circ t_n : t_n^{-1} (X_0) \to
\tau_n^{-1} (Y_0)$.

Since $t_n^{-1} (X_0)$ are isometric copies of $X_0$ containing
$x_0$, by Ascoli Theorem they converge to an isometric copy $X_1$
of $X_0$. A similar argument can be done for $\tau_n^{-1} (Y_0)$,
which converge to an isometric copy $Y_1$ of $Y_0$, therefore
$\q_n$ converges to an $(L,C)$--quasi-isometry $\tilde{\q} : X_1
\to Y_1$. Also, the extensions $\q_n^e =\tau_n^{-1} \circ \q_e
\circ t_n$ of $\q_n$ to $\hip_\R^3$ converge to an extension
$\tilde{\q}_e$ of $\tilde{\q}$.

On the other hand, since $t_n$ and $\tau_n$ restricted to $\C
\subset \partial_\infty \hip^3_\R $ are homotheties of center $O$,
the restrictions of the boundary maps $\partial \q_n^e : \C \to \C
$ converge to the differential $d_O \partial \q_e$. Thus, $d_O
\partial \q_e$ is the restriction to $\C$ of $\partial \tilde{\q}_e$. From
this it can be deduced that $d_O \partial \q_e$ is a similarity.
We give the sketch of proof below. The full proof is more
elaborate and can be found in \cite{Sch}.

The argument in Step 1 implies that $\tilde{\q}$ sends every
boundary horosphere of $X_1$ at uniformly bounded distance of a
boundary horosphere of $Y_1$. From this it can be deduced that all
horospheres having a certain Euclidean height $h$ are sent at
uniformly bounded distance from horospheres having an Euclidean
height in $[h/\lambda , \lambda h ]$ (for some constant $\lambda
\geq 1$ depending on the constant $M$ given by Step 1 for
$\tilde{\q}$). Note that the basepoints of horospheres of heights
at least $h$ in $X_1$ compose nets $\mathbf{N}_h$ in $\C$ with the
corresponding constants $\delta$ and $\epsilon$ smaller and
smaller as $h$ decreases to zero. Thus, $d_O
\partial \q_e$ sends each of these nets $\mathbf{N}_h$ of $\C$ into
a net $\mathbf{N}_{h/\lambda}$ of $\C$. Up to now, nothing
surprising, since $d_O \partial \q_e$ is a linear map.

Now we change stereographic projection, and in both the domain and
the range model reverse $O$ with $\infty$. In the new models, we
again have $\tilde{\q} :X_1 \to Y_1 $ extended to $\tilde{\q}_e
:\hip_\R^3 \to \hip_\R^3 $ such that $\partial \tilde{\q}_e$ fixes
both $O$ and $\infty$ and such that its restriction to $\C
\setminus \{
 O\}$ coincides with $I\circ d_O \partial \q_e \circ I$, where $I$
 is the inversion with respect to the unit circle. An
 argument as above implies that $I\circ d_O \partial \q_e \circ I$
 sends nets $\mathbf{N}_h$ of $\C $ into nets $\mathbf{N}_{h/\lambda}$ of $\C$ for every $h$.
  This implies that $I\circ d_O \partial \q_e \circ I$ is also
  linear. But this can happen only if $d_O \partial \q_e$ is a
  similarity.

  \me

  \textsc{Step 4.} Theorem \ref{mob} and Step 3 imply that there
  exists an isometry $g$ of $\hip^3_\R$ such that $\partial g = \partial
  \q_e$. It follows that $g$ and $\q_e$ are at uniformly bounded
  distance from each other. In particular $g$ restricted to $X_0$
  is at uniformly bounded distance from $\q$. Next it is shown that $g$ is
  in the commensurator $Comm
(\Gamma ,G)$ of $\Gamma $ into $G$.

The argument is by contradiction. Suppose that $\Gamma^g \cap G$
has infinite index either in $\Gamma^g$ or in $ G$. Without loss
of generality we may assume that it has infinite index in $G$. It
follows that there exists a sequence $g_n$ of elements in $G$ such
that $g_n (\Gamma^g)$ are distinct left cosets in the group of
isometries of $\hip^n_\F$.

Let $\beta $ be a basepoint of a boundary horosphere of $Y_0$. It
is an easy exercise to show that, up to taking a subsequence,
there exists a sequence $\gamma_n$ in $\Gamma $ and another
basepoint $\alpha $ of a boundary horosphere of $Y_0$ such that
$g_n \gamma_n^g (\alpha )=\beta $ and $g_n \gamma_n^g (y_0)\in
B(y_0, R )$, where $R$ is a constant.

Consider the respective stereographic projections sending $(\alpha
, y_0)$ to $(\infty , (0,0,1))$ on the $\hip^3_\R$ of definition
and $(\beta , y_0)$ to $(\infty , (0,0,1))$ on the range
$\hip^3_\R$. In these new half-space models of $\hip^3_\R$ the
isometry $g_n \gamma_n^g$ fixes $\infty$ and sends $(0,0,1)$ at
distance at most $R$ from itself. Also, since $\gamma_n^g$ are
isometries of $g(X_0)$, $g_n$ are isometries of $Y_0$ and $Y_0$ is
at uniformly bounded distance from $g(X_0)$, it follows that $g_n
\gamma_n^g (Y_0)$ is at Hausdorff distance at most $D$ from $Y_0$,
where $D$ is a constant independent of $n$.

By Ascoli Theorem, $g_n \gamma_n^g $ converges to an isometry
$\hat{g}$ such that $Y_1=\hat{g}(Y_0)$ is at distance at most $D$
from $Y_0$.

Let $\mathbf{N}_h$ be the set of basepoints of boundary
horospheres of $Y_0$ of Euclidean height at least $h$. Note that
$G_\infty$, the stabilizer in $G$ of the point $\infty$, acts on
$\C$ such that $G_\infty \backslash \C$ is a flat torus. Let $D
\subset \C$ be a fundamental domain (quadrangle) projecting on
this torus. The number of horoballs of $Y_0$ of Euclidean height
at least $h$ and with basepoints in $D$ is finite. Let
$\mathbf{B}_h$ be the finite set of their basepoints.

Then $\mathbf{N}_h =\coprod_{b\in \mathbf{B}_h} G_\infty b$ is a
finite union of grids in $\C$, hence a net in $\C$.

The previous considerations imply that $\mathbf{N}_h^n=g_n
\gamma_n^g (\mathbf{N}_h)$ is a net contained in
$\mathbf{N}_{h/\lambda }$ and likewise for the net
$\hat{\mathbf{N}}_h=\hat{g} (\mathbf{N}_h)$. On the other hand
$\mathbf{N}_h^n$ converges to $\hat{\mathbf{N}}_h$ in the
compact-open topology. The only way in which this convergence can
occur, both nets being in the larger net $\mathbf{N}_{h/\lambda
}$, is that they coincide on larger and larger subsets.

An isometry of $\hip^3_\R $ fixing four points on $\partial_\infty
\hip^3_\R $ which are not on the same circle is the identity
isometry\fn{Elementary proof: an isometry fixing three distinct
points in the boundary has to fix a point in $\hip^3_\R$
\cite[Proposition A.5.14]{BP}. Thus the isometry can be identified
with a matrix in $SO(3)$ fixing four tangent vectors-the vectors
pointing towards the four fixed points in $\partial_\infty
\hip^3_\R $. If the four vectors were in the same plane then the
corresponding points in the boundary would be on the same circle.
Therefore three of the four vectors are linearly independent and
the isometry has to be the identity.}. Hence two isometries which
coincide on four points on $\partial_\infty \hip^3_\R $ not on the
same circle are equal. It follows that the sequence $g_n
\gamma_n^g$ becomes stationary, for $n$ large enough. This
contradicts the hypothesis that $g_n (\Gamma^g)$ are distinct left
cosets.\hspace*{\fill}$\Box$

\subsection{Proof of Theorem \ref{RSch}}

A non-uniform lattice of isometries of $\hip^n_\F$ has infinitely
many ends if and only if $\hip^n_\F=\hip^2_\R$. On the other hand,
having infinitely many ends is a quasi-isometry invariant. Thus,
either both $X_1$ and $X_2$ are $\hip^2_\R$ or both differ from
it. Suppose we are in the second case. The quasi-isometry $\q :G_1
\to G_2$ induces as in the previous section a quasi-isometry $\q
:X_1 \to X_2$, hence a quasi-conformal homeomorphism between the
boundaries at infinity $\partial_\infty X_1$ and $\partial_\infty
X_2$. Elementary dimension and structure arguments imply that
$X_1=X_2$.

The rest of the statement (1) follows from the Theorem \ref{rig}.

We now get to the proof of statement (2). The following standard
fact is needed.

\begin{lemma}\label{quasia}
Let $\Lambda $ and $G$ be finitely generated groups and let $\q
:\Lambda \to G$, $\bar{\q }: G \to \Lambda$ be two quasi-converse
$(L_0,C_0)$--quasi-isometries. Then to every $\lambda \in \Lambda
$ one can associate an $(L,C)$--quasi-isometry of $G$, $\q_\lambda
= \q\circ L_\lambda \circ \bar{\q}$, where $L_\lambda$ denotes the
isometry on $\Lambda $ determined by the action of $\lambda $ to
the left, and $(L,C)$ can be effectively computed from
$(L_0,C_0)$. Moreover the map $\lambda \to \q_\lambda $ defines a
\begin{itemize}
    \item \textbf{quasi-action of $\Lambda $ on $G$}: there exists $D=D(L_0,C_0)$ so that
     for every $\lambda , \eta \in \Lambda $ the following holds:
    \begin{eqnarray}
      \dist (\q_\lambda \circ \q_\eta \, ,\, \q_{\lambda \eta }) &\leq & D\, , \\
      \dist (\q_\lambda \circ \q_{\lambda^{-1}} \, ,\, \id ) &\leq &
      D\, .
    \end{eqnarray}
    \item \textbf{which moreover is quasi-transitive}: for
    every $g,g'\in G$ there exists $\lambda \in \Lambda $ such
    that
    $$
\dist (\q_\lambda (g), g')\leq C_1\, ,
    $$ where $C_1=C_1(L_0,C_0)$;
    \item \textbf{and of finite quasi-kernel}: for every
    $K>0$ the set of $\lambda $ such that $\dist (\q_\lambda , \id)\leq
    K$ is finite.
\end{itemize}
\end{lemma}

The proof of the lemma is left as an exercise to the reader.

\me

In the particular case considered, Lemma \ref{quasia} and Theorem
\ref{rig} imply that if $\Lambda $ is quasi-isometric to $G$
non-uniform lattice of isometries of $X=\hip^n_\F\neq \hip^2_\R$,
then there exists a morphism of finite kernel
$$
\phi :\Lambda \to Comm (G)\, .
$$

The fact that the image $G_1$ of $\Lambda $ under $\phi$ is
discrete can be proved by contradiction. Suppose it is not
discrete, hence there are infinitely many elements in $\phi
(\Lambda )$ in the neighborhood of the identity element $\id\in
G$. Then for some $D$ large enough we have that for infinitely
many $\lambda \in \Lambda$
$$
\dist (\q \circ L_\lambda \circ \bar{\q }(\id) , \id )\leq D
\Rightarrow \dist (L_\lambda \circ \bar{\q }(\id) , \bar{\q }(\id)
)\leq D'\, ,
$$ where $D'=D'(L,C,D)$. This contradicts the
fact that every ball in $\Lambda$ is a finite set.

Also, one can argue that $G_1\backslash X$ has finite volume
roughly as follows. Consider a complementary set $X_0$ in $X$ of a
family of countably many pairwise disjoint open horoballs such
that $G\backslash X_0$ is compact. Lemma \ref{quasia} implies that
$\Lambda$ acts quasi-transitively by quasi-isometries on $G$,
hence on $X_0$. It follows that $G_1$ acts ``with compact
quotient'' on $X_0$. See \cite[$\S 10.4$]{Sch} for details.
\endproof

\section{Classes of groups complete with respect to
quasi-isometries}\label{cgc}

Another way of interpreting Theorem \ref{RSch}, (2), is the
following. Let $\mathcal{C}$ be the class of non-uniform lattices
of isometries in $\hip^n_\F \neq \hip_\R^2$. Then every finitely
generated group $\Lambda$ quasi-isometric to a group $G\in
\mathcal{C}$ is itself in $\mathcal{C}$, up to taking its quotient
by a finite normal subgroup. One may ask what other classes of
groups behave similarly. Possibly, to the operation of taking
quotient by finite normal subgroup one has to add the other
algebraic operation preserving the quasi-isometry class: taking a
subgroup of finite index.

\begin{definition}\label{qic}
A class of finitely generated groups $\mathcal{C}$ with the
property that if $\Lambda $ is quasi-isometric to $G\in
\mathcal{C}$ then $\Lambda_1\in \mathcal{C}$, where $\Lambda_1$ is
either a finite index subgroup of $\Lambda $ or a quotient of
$\Lambda $ by a finite normal subgroup, is called \textit{class of
groups complete with respect to quasi-isometries} or \textit{q.i.
complete}.
\end{definition}

The question of finding such classes has been asked for the first
time by M. Gromov in \cite{Gr1}.

\subsection{List of classes of groups q.i.
complete}

We give a (non-exhaustive) list of classes of groups q.i.
complete. We begin with the q.i. complete classes of lattices of
isometries of symmetric spaces other than those discussed above.
All the lattices that we consider are supposed to be
\textit{irreducible}.

\begin{enumerate}
    \item uniform lattices of isometries of a symmetric space $X$
    for the list of spaces $X$ below.
      \begin{itemize}
        \item $X=\hip_\R^n$, with $n\geq 3$, by the work of Sullivan
and Tukia (see the lecture notes of Marc Bourdon and references
therein);
         \item $X=\hip^n_\hip$ and $X= \hip^2_{Cay}$ by the work
         of P. Pansu
    \cite{Pa2};
    \item $X=\hip^n_\C\, ,n\geq
    2$, by the work of R. Chow \cite{Cho};
    \item $X=\hip_\R^2$. In this case a proof of the q.i.
    completeness goes as follows:
    \begin{itemize}
        \item A finitely generated group $\Lambda $
        quasi-isometric to $\hip^2_\R $ is a hyperbolic group (see
        Example 6 below for a definition), with boundary at
        infinity homeomorphic to $\sph^1$. Every
        hyperbolic group acts on its boundary as a convergence
        group \cite{Tu2}. For a definition of convergence groups
        see Section \ref{dic}.
        \item Every convergence group is conjugate to a Fuchsian group in
        Homeo$(\sph^1)$. This follows from \cite{Tu1}, \cite{CJ}
        and \cite{Ga}.
    \end{itemize}
        \item $X$ irreducible symmetric
space of rank at least $2$. This result is due to B. Kleiner and
B. Leeb \cite{KlL}. See also \cite{EF} for another proof.
      \end{itemize}
    \item non-uniform lattices of isometries of a symmetric space $X$ of rank at least
    $2$. This is due to R. Schwartz for a family of $\Q$--rank one
    lattices containing the Hilbert modular groups
    \cite{Sch2} and to A. Eskin \cite{E} in the general case, under the condition that
     $X$ has no factors of rank $1$. See
    also \cite{Dr2} for another proof of the general case.

    Moreover, in this case Statement (1), (b), of Theorem
    \ref{RSch} holds as well, that is: two non-uniform lattices are
    quasi-isometric if and only if they are commensurable.

\begin{remark}
In the cases of uniform lattices in $X=\hip^n_\hip$ and $X=
\hip^2_{Cay}$ or $X$ of rank at least two, as well as in the case
of non-uniform lattices of a symmetric space $X$ of rank at least
$2$, the q.i. completeness result is obtained via a rigidity
result similar to Theorem \ref{rig}, that is : a quasi-isometry of
the lattice is at bounded distance from an isometry. Moreover, in
the case of non-uniform lattices, this isometry is in the
commensurator of the lattice.
\end{remark}
    \item fundamental groups of non-geometric Haken manifolds with zero Euler characteristic \cite{KaL2}
    (see the lecture notes of M. Kapovich).
    \item finitely presented groups \cite{GH0}.
    \item nilpotent groups. This follows from the Polynomial Growth Theorem of M. Gromov
    \cite{Gr3}. We recall that the growth function  $B_S:\N \to \N$ of a group $G$
with a finite set of generators $S$ is defined by $B_S(n)=$ the
cardinal of the ball $B_S(1,n)$ in the word metric $\dist_S$. The
theorem of M. Gromov states that the growth function with respect
to some (hence any) finite set of generators is polynomial if and
only if the group is virtually nilpotent.

    The subclass of Abelian groups is also q.i. complete, as follows from results in
    \cite{Pa}. See the discussion in Section \ref{exac}.

    \item hyperbolic groups.

    We recall that a geodesic metric space is called
    $\delta$-\textit{hyperbolic} if in every geodesic triangle, each edge
    is contained in the $\delta$-tubular neighborhood of the
    union of the other two edges. If $\delta=0$ the space is
    called \textit{real tree} or \textit{$\R$--tree}.

    A finitely generated group is called \textit{hyperbolic} if
    its Cayley graph is hyperbolic. For instance, uniform lattices
    in rank one symmetric spaces are such.

    The q.i. completeness of the class of hyperbolic groups
    follows easily from the definition and from

    \begin{lemma}[Morse lemma, see \cite{GH}]
Every $(L,C)$--quasi-isometric segment in a $\delta$-hyperbolic
space is at Hausdorff distance at most $D$ from the geodesic
segment joining its endpoints, where $D=D(L,C,\delta)$.
    \end{lemma}
    \item amenable groups \cite{GH0}. We recall that a discrete group $G$ is \textit{amenable}
if for every finite subset $K$ of $G$ and every $\epsilon \in
 (0,1)$ there exists a finite subset $F\subset G$ satisfying:
$$
\mathrm{card}\, KF < (1+\epsilon ) \mathrm{card}\, F\, .
$$

    \item the whole class of solvable groups is not q.i. complete, as pointed out
by the counter-example in \cite{Dyu}. Still, there are some
smaller classes of solvable groups that are q.i. complete. See for
instance \cite{FM}, \cite{FM2}, \cite{EFW1},  \cite{EFW2}.
\end{enumerate}

\subsection{Relatively hyperbolic groups: preliminaries}

In the same way in which uniform lattices in rank one symmetric
spaces inspired the notion of hyperbolic group, non-uniform
lattices inspired the notion of \textit{relatively hyperbolic
group}.\fn{What is called in this paper relatively hyperbolic
group is sometimes called in the literature \textit{strongly
relatively hyperbolic group}, in contrast with \textit{weakly
relatively hyperbolic group}.} This notion was defined by M.
Gromov in \cite{Gr2}. Then several equivalent definitions of the
same notion as well as developments of the theory of relatively
hyperbolic groups were provided in \cite{Bowrh}, \cite{Fa},
\cite{Dahmani:thesis}, \cite{Yaman:RelHyp}, \cite{DS},
\cite{Osin}, \cite{Dr3}.

Here we recall the definition of B. Farb \cite{Fa}. Let $G$ be a
finitely generated group and let $\{ H_1,\dots , H_m \}$ be a
collection of subgroups of $G$. Let $S$ be a finite generating set
of $G$ invariant with respect to inversion. The idea is to write
down a list of properties which force $G$ to behave with respect
to $\{ H_1,\dots , H_m \}$ in the same way in which a non-uniform
lattice $\Gamma $ behaves with respect to its cusp subgroups $\{
\Gamma_{\alpha_1},\dots , \Gamma_{\alpha_m} \}$ (see Figure
\ref{fig4}).

We denote by $\mh$ the set $\bigsqcup_{i=1}^m (H_i \setminus \{ 1
\})$. We can consider two Cayley graphs for the group $G$, $\cgs$
and $\cgsh$. We note that $\cgs$ is a subgraph of $\cgsh$, with
the same set of vertices but a smaller set of edges, and that
$\cgsh$ is not locally finite if at least one of the subgroups
$H_i$ is infinite. We have that $\dsh(u,v)\leq \ds (u,v)$, for
every two vertices $u,v$.

\me

\begin{definition}
Let $\pgot$ be a path in $\cgsh$. An $\mh$--\textit{component} of
$\pgot$ is a maximal sub-path of $\pgot$ contained in a left coset
$gH_i,\, i\in \{1,2, \dots ,m \},\, g\in G$.

The path $\pgot$ is said to be \textit{without backtracking} if it
does not have two distinct $\mh$--components in the same left
coset $gH_i$.
\end{definition}

The notion of weak relative hyperbolicity has been introduced by
B. Farb in \cite{Fa}. We use a slightly different but equivalent
definition. The proof of the equivalence can be found in
\cite{Osin}.

\begin{definition}
The group $G$ is \textit{weakly hyperbolic relative to} $\{
H_1,\dots , H_m \}$ if and only if the graph $\cgsh$ is
hyperbolic.
\end{definition}

This property is not enough to determine a picture as in Figure
\ref{fig4}. For instance $G=\Z^2$ satisfies the previous property
with respect to its subgroup $H=\Z \times \{ 0\}$. This case does
not at all look as in Figure \ref{fig4}, in that the tubular
neighborhoods of left cosets of $H$ do not, as in Figure
\ref{fig4}, intersect in a finite set, but on the contrary the
intersection may contain both left cosets. Vaguely speaking, in
Figure \ref{fig4} the left cosets stay close in the respective
neighborhoods of a pair of points realizing the minimal distance,
and then diverge, while in the example above two left coset stay
parallel.

One has to add a second property in order to obtain the proper
image, and thus to define (strong) relative hyperbolicity. Before
formulating this property, we must mention another notable example
of group weakly relatively hyperbolic and not strongly relatively
hyperbolic. The Mapping Class Group of a hyperbolic surface
$\Sigma$ (also defined as $Out (\pi_1 (\Sigma))$) is weakly
hyperbolic (and \textit{not} strongly hyperbolic) relative to
finitely many stabilizers of closed geodesics on the surface
$\Sigma$. The weak relative hyperbolicity follows from \cite{MM}
(see also \cite{Bow}). The reason for which (strong) relative
hyperbolicity is not satisfied is again that the intersection of
two tubular neighborhoods of left cosets is not finite. Indeed,
two stabilizers of two closed geodesics intersect in the
stabilizer of both, which can be itself infinite.

\me

\Notat \quad For every path $\pgot$ in a metric space $X$, we
denote the start of $\pgot$ by $\pgot_-$ and the end of $\pgot $
by $\pgot_+$. \me

\me

\begin{definition}
The pair $(G\, ,\, \{ H_1,\dots , H_m \})$ satisfies the
\textit{Bounded Coset Penetration (BCP) property} if for every
$\lambda \geq 1$ there exists $a=a(\lambda )$ such that the
following holds. Let $\pgot$ and $\q$ be two
$\lambda$-bi-Lipschitz paths without backtracking in $\cgsh$ such
that $\pgot_-=\q_-$ and $\ds (\pgot_+,\q_+)\leq 1$.
\begin{itemize}
  \item[(1)] Suppose that $s$ is an $\mh$--component of $\pgot$ such that $\ds (s_-,s_+)\geq
  a$. Then $\q$ has an $\mh$--component contained in the same left
  coset as $s$;
  \item[(2)] Suppose that $s$ and $t$ are two $\mh$--components of $\pgot$ and $\q$, respectively,
   contained in the same left
  coset. Then $\ds (s_-,t_-)\leq a$ and $\ds (s_+,t_+)\leq a$.
\end{itemize}
\end{definition}

In particular BCP property implies that if $H_i$ is infinite two
left cosets of $H_i$ cannot be at finite Hausdorff distance one
from the other. The proof is left as an exercise to the reader.

\begin{definition}
The group $G$ is \textit{(strongly) hyperbolic relative to} $\{
H_1,\dots , H_m \}$ if it is weakly hyperbolic relative to $\{
H_1,\dots , H_m \}$ and if $(G\, ,\, \{ H_1,\dots , H_m \})$
satisfies the BCP property.
\end{definition}

Both in the case of weak and strong relative hyperbolicity, the
subgroups $H_1,\dots , H_m$ are called \textit{peripheral
subgroups}. A subgroup conjugate to some $H_i,\, i\in
\{1,...,m\}$, is called a \textit{maximal parabolic subgroup}. A
subgroup contained in a maximal parabolic subgroup is called
\textit{parabolic}.

\me

\textit{Other examples of relatively hyperbolic groups (besides
non-uniform lattices)}:
\begin{enumerate}
    \item $A*_F B$, where $F$ is finite, is hyperbolic relative to $A$ and $B$; more
    generally, fundamental groups of finite graphs of groups with
    finite edge groups are hyperbolic relative to the vertex
    groups \cite{Bowrh}.
    \item a hyperbolic group $\Gamma $ is hyperbolic relative to
     \begin{itemize}
        \item $H=\{ 1\}$;
        \item any class of infinite quasi-convex subgroups $\{ H_1,\dots , H_k \}$
        with the property that $H_i^g \cap H_j$ is finite if $i\neq j$ or $g\not\in
        H_i$ \cite[Theorem 7.11]{Bowrh}.

        For instance let $\Gamma$ be a uniform lattice of
        isometries of $\hip^3_\R$ such that for some totally geodesic
        copy of the hyperbolic plane $\hip^2_\R$ in $\hip^3_\R$,
        $H=\Gamma \cap Isom (\hip^2_\R)$ is a uniform lattice of
        $\hip^2_\R$. Then $\{ H \}$ satisfies the previous
        properties.
     \end{itemize}
    \item fundamental groups of complete finite volume manifolds of pinched
    negative sectional curvature are hyperbolic relative to the
    fundamental groups of their cusps (\cite{Bowrh}, \cite{Fa});
    \item fundamental groups of (non-geometric) Haken manifolds with at least one
    hyperbolic component are hyperbolic relative to fundamental
    groups of maximal graph-manifold components and to fundamental groups of
    tori and Klein bottles not contained in a graph-manifold
    component; this follows from the previous example and from the
    combination theorem of F. Dahmani \cite[Theorem 0.1]{Dah} (for a
     combination theorem applying also to non-finitely generated
     groups see \cite{Osin:combination});
    \item  fully residually free groups, also known as limit
    groups, are hyperbolic relative to their maximal
    Abelian non-cyclic subgroups \cite{Dah}. Moreover, according to \cite{AB} these
    groups are known to be $CAT(0)$ with isolated flats, after the terminology in \cite{Hr};
    \item more generally, finitely generated groups acting freely
    on $\R^n$--trees are hyperbolic relative to their maximal
    Abelian non-cyclic subgroups \cite{Gui}.
\end{enumerate}

\me

\begin{remark}
Throughout the discussion of relatively hyperbolic groups we
tacitly rule out the case of a finite group hyperbolic relative to
any class of subgroups, as well as the case when one of the
peripheral subgroups $H_i$ is the ambient group itself.

Thus, we are in the case of an \emph{infinite} group $G$ and a
finite collection (possibly reduced to one element) $\{ H_1,\dots
,H_m\}$ of \emph{proper} subgroups of $G$. In this case it follows
that each subgroup $H_i$ has infinite index in $G$. We note also
that if all peripheral subgroups are finite then $G$ is
hyperbolic; if there exists at least one infinite peripheral
subgroup, then the finite peripheral subgroups can be removed from
the list, and it can be thus assumed that all $H_i$ are infinite.
\end{remark}

\begin{remark}
Recently, relatively hyperbolic groups have been used to construct
examples of infinite groups with exotic properties. Thus in
\cite{Osin:twoconjugacy} it is proved that there exist uncountably
many pairwise non-isomorphic two-generated groups without finite
subgroups and with exactly two conjugacy classes, answering an old
question in group theory.
\end{remark}

\me

\begin{ques}\label{qrigrh}
Is the class of relatively hyperbolic groups q.i. complete ?
\end{ques}

\me

Before discussing Question \ref{qrigrh}, we define our main tools.

\section{Asymptotic cones of a metric space}

\subsection{Definition, preliminaries}\label{prel}

The notion of asymptotic cone was defined in an informal way in
\cite{Gr3}, and then rigorously in \cite{VDW} and \cite{Gr4}. The
idea is to construct, for a given metric space, an image of it
seen from infinitely far away.

First one needs the notion of \textit{non-principal ultrafilter}.
This can be defined as a finitely additive measure $\omega$
defined on the set of all subsets of $\N$ (or, more generally, of
a countable set) and taking only values zero and one, such that on
all finite subsets it takes value zero. In particular if $\N = A_1
\sqcup \cdots \sqcup A_n$ and all $A_i$ are infinite, then there
exists $i_0\in \{ 1,2,\dots , n\}$ such that $\omega (A_{i_0})=1$
and $\omega (A_{j})=0$ for every $j\neq i_0$.

The fact that $\omega $ takes only values $0$ and $1$ immediately
brings to one's mind the idea of a characteristic function.
Indeed, $\omega$ satisfies the previous properties if and only if
it is the characteristic function $\unu_{\mathcal{U}}$ of a
collection $\mathcal{U}$ of subsets of $\N $ which is
\begin{itemize}
    \item an \textit{ultrafilter}, that is, a
maximal filter;
    \item moreover \textit{nonprincipal}, that is containing \textit{the Fr\'echet filter}.
\end{itemize}

For definitions of the notions above, that is, for a list of
axioms, see Section \ref{dic}. For more details see \cite{Bou}.
The main advantage of the second way of defining non-principal
ultrafilters resides, besides the questionable pleasure of dealing
with axioms, in the fact that it shows that such objects always
exist. We also note that a functional analytic treatment of
ultrafilters is possible. Thus, the notion of $\omega$--limit can
be seen as an application of the Hahn-Banach theorem to the space
of relatively compact sequences, the subspace of convergent
sequences and the limit map on it.

Since all ultrafilters in this paper are nonprincipal, we drop
this adjective henceforth.

Given an ultrafilter $\omega$ and a sequence $(x_n)$ in a
topological space, one can define the $\omega$--\textit{limit}
$\lim_\omega x_n$ of the sequence as an element $x$ such that for
every neighborhood $\nn$ of it, $$\omega \left( \{ n\in \N \mid
x_n\in \nn \} \right)=1\, .$$

The following property emphasizes the main interest of
(nonprincipal) ultrafilters.

\begin{proposition}\cite{Bou}
If $(x_n)$ is contained in a compact, its $\omega$--limit always
exists.
\end{proposition}

Note that, as soon as it exists, the $\omega$--limit is also
unique.
 Also, it is not difficult to see that it is a limit of a
 converging subsequence. Thus, an ultrafilter is a device to
 select a point of accumulation for any relatively compact sequence, in a coherent manner.
 In some sense, it is a systematic approach to the process of taking the diagonal subsequence,
  after selecting converging subsequences in countably many sequences.

With such a tool at hand, which makes almost any reasonable
sequence converge, one can hope to define, for a given metric
space $(X,\dist )$, an image of it seen from infinitely far away.
More precisely, one has to take a sequence of positive numbers
$d_n$ diverging to infinity, and try to construct a limit of the
sequence of metric spaces $\left( X , \frac{1}{d_n}\dist \right)$.

As in the formal construction of completion, one can simply take
the set $\mathcal{S}(X)$ of all sequences $(x_n)$ in $X$ and try
to define a metric on this space by
$$
\dist_\omega \left( x,y \right)=\lim_\omega \frac{\dist
(x_n,y_n)}{d_n},\, \mbox{ for }x=(x_n),\, y=(y_n)\, .
$$

The problem is that the latter limit can be $+\infty $, or it can
be zero for two distinct sequences.

To avoid the situation $\dist_\omega \left( x,y \right)=+\infty $,
one restricts to a subset of sequences defined as follows. For a
fixed sequence $e=(e_n)$, consider
\begin{equation}\label{se}
 \mathcal{S}_e(X) = \left\{ (x_n) \in X^\N \: ;\:  \left(
\frac{\dist (x_n, e_n)}{d_n}\right) \mbox{ is a bounded sequence}
\right\}\, .
\end{equation}

To deal with the situation when $\dist_\omega \left( x,y
\right)=0$ while $x\neq y$, one uses the classical trick of taking
the quotient for the equivalence relation
$$
x\sim y \Leftrightarrow \dist_\omega \left( x,y \right) =0\, .
$$

The quotient space $\mathcal{S}_e (X) / \sim $ is denoted $\co{X;
e,d }$ and it is called \textit{the asymptotic cone of $X$ with
respect to the ultrafilter $\omega$, the scaling sequence
$d=(d_n)$ and the sequence of observation centers $e$}.

A sequence of subsets $(A_n)$ in $X$ gives rise to a \textit{limit
subset} in the cone, defined by
$$
\lio{A_n} = \{ \lio{a_n} \mid a_n \in A_n,\,  \forall \, n\in \N
\}\, .
$$

If $\lim_\omega \frac{\dist (e_n,A_n)}{d_n}=+\infty $ then
$\lio{A_n}=\emptyset$.

\bigskip

\noindent\textit{Properties of asymptotic cones:}

\begin{enumerate}
    \item $\co{X; e,d}$ is a complete metric space;
    \item every limit subset $\lio{A_n}$, if non-empty, is closed;
    \item if $X$ is geodesic then every asymptotic cone $\co{X; e,d}$
    is geodesic;
    \item an $(L,C)$--quasi-isometry between two metric spaces $\q : X\to
    Y$ gives rise to a bi-Lipschitz map between asymptotic cones

\begin{eqnarray*}
  \q_\omega : \co{X;e,d} & \to & \co{Y;\q(e),d} \\
  \lio{x_n} &\to & \lio{\q (x_n)}\: ;
\end{eqnarray*}

    \item If $G$ is a group then every $\co{G;e,d}$ is isometric
    to $\co{G; 1,d}$, where $1$ denotes here the constant sequence
    equal to $1$;
    \item $\co{G; 1,d}$ is a homogeneous space.
\end{enumerate}

Proofs of the previous properties can be found in \cite{Gr4},
\cite{KlL}, \cite{KaL1}. None of them is difficult; they are good
exercises in order to get familiar with the notion. We shall take
a closer look only at the last property, namely we shall exhibit
the group acting transitively by isometries on $\co{G; 1,d}$. Let
$G^\N$ be the set of all sequences in $G$ and let
$\mathcal{S}_{(1)}(G)$ be the subset of sequences defined as in
(\ref{se}). We consider the equivalence relation
$$
(g_n)\approx (g_n')\Leftrightarrow \omega (\{ n\in \N \mid
g_n=g_n'\})=1\, .
$$

The quotient space $\Pi_\omega G = G^\N / \approx $ is a group,
called the $\omega$--\textit{ultrapower of $G$}. The subgroup
$G^\omega = \mathcal{S}_{(1)}(G)/\approx $ acts transitively by
isometries on $\co{G;1,d}$ by:
$$(g_n)^\omega \lio{x_n}=\lio{g_nx_n}\, .$$

\me

One can put a condition in order to restrict the growth of the
scaling sequence with respect to the ultrafilter. The idea is to
choose a sequence $d$ and an ultrafilter $\omega$ such that there
is no set $E$ with $\omega (E)=1$ and such that $(d_n)_{n\in E}$
grows faster than exponentially. The definition is as follows:

\begin{definition}
The pair $(\omega , d)$ is \textit{non-sparse} if:
\begin{itemize}
    \item[(i)] For every $a>1$ we have $d_n\leq a^n$
    $\omega$---almost surely.
    \item[(ii)] For every ordered infinite subset $E=\{ i_1,i_2,\dots , i_n,\dots
    \}$ such that there exists $a>1$ satisfying $\lim_{n\to \infty }
    \frac{d_{i_n}}{a^n}=+\infty$,
    we have that $\omega (E)=0$.
\end{itemize}

In one of the properties (i) and (ii) is not satisfied, we say
that the pair $(\omega , d)$ is \textit{sparse}.
\end{definition}

Both sparse and non-sparse pairs exist. In order to construct a
sparse pair, it suffices to take an ultrafilter $\omega$, and via
the injection $n\mapsto 3^n$ to identify it to an ultrafilter
supported by the set $\{ 3^n \; ;\; n\in \N \}$.

To construct a non-sparse pair, take for instance $d_n=n$ and all
the sets $E$ described in (ii), for this choice. The collection of
subsets of $\N$ having complementaries in $\N$ either finite or
contained in a set of type $E$ is a filter. An ultrafilter
$\mathcal{U}$ containing it is non-principal and the pair
$\left(\, \unu_{\mathcal{U}}\, ,\, (n)\, \right)$ is non-sparse.

\subsection{A sample of what one can do using asymptotic cones}

\begin{proposition}[\cite{Gr4}]\label{sample}
Let $\Gamma$ be a discrete group endowed with a metric $\dist $,
left invariant with respect to the action of the group on itself,
     such that all balls are finite.
     \begin{itemize}
        \item[(1)] If all the asymptotic cones of $(\Gamma , \dist )$ are path-connected then $\Gamma $ is finitely
     generated.
        \item[(2)] If moreover all the asymptotic cones of $(\Gamma , \dist )$ are simply connected then $\Gamma $ is finitely
     presented.
     \end{itemize}
\end{proposition}

\begin{remark}
\begin{itemize}
    \item[(a)] A metric as in Proposition \ref{sample}, (1), can be obtained for instance if the
     group acts properly discontinuously and freely by isometries on a
     proper metric space $X$. Given $x\in X$ we identify $\Gamma $
     with the orbit $\Gamma x$ and we take the induced metric.

     A particular case of the previous situation is when
     $\Gamma$ is a subgroup of a finitely generated group $G$.
     Then we can take as $X$ the Cayley graph of $G$. If we choose $x=1$, the induced metric $\dist $
     is the word metric of $G$ restricted to $\Gamma$.
    \item[(b)] The converse of Proposition \ref{sample}, (2), is
    not true. This can be seen for instance in the case of
    Baumslag-Solitar groups $BS(p,q)$ or in the case of uniform
    lattices in the solvable group $Sol$. These groups are finitely
    presented, nevertheless their asymptotic cones have
    uncountable fundamental group \cite{Bu}.
\end{itemize}
\end{remark}

     \proof \textbf{(1)}\quad \textsc{Step} 1. We first prove that between every
pair of elements $x,y $ in $\Gamma
     $ there exists a discrete path composed of at most $N$ steps of length at
     most $\frac{\dist (x,y)}{M}$, where both $M$ and $N$ are fixed.

     Here is the precise statement with quantifiers: for every $M>1$ there
     exist $N\in \N ,\, N\geq 2,$ and $D>0$ such that for every $x,y\in \Gamma $
     with $\dist (x,y)\geq D$, there exists a finite sequence of
     points $t_0=x,t_1,\dots , t_m=y$ with $m\leq N$ and $\dist (t_i,t_{i+1})\leq \frac{\dist
     (x,y)}{M}$ for every $i\in \{ 0,1,\dots , m-1\}$.

     We argue by contradiction and suppose it is not the case. Then there exists $M$ and a sequence of pairs of points $(x_n,y_n)\in \Gamma \times \Gamma
     $ with $d_n=\dist (x_n,y_n )\geq n$ and such that every discrete
     path of at most $n$ steps between $x_n$ and $y_n$ has at least one step of length $\geq \frac{\dist
     (x_n,y_n)}{M}$. In the asymptotic cone $\co{\Gamma ; x_n , d_n
     }$ the two sequences $(x_n)$ and $(y_n)$ give two points
     $x_\omega$ and $y_\omega$ at distance 1 such that for every
     $n$, every discrete path joining  $x_\omega$ and $y_\omega$ and having $n$ steps
     has at least one step of length $\geq \frac{1}{M}$. On
     the other hand, since $\co{\Gamma ; x_n , d_n
     }$ is path-connected, $x_\omega$ and $y_\omega$ can be joined
     by a path. On this path can be chosen a finite discrete path of
     steps at most $\frac{1}{M}$ between $x_\omega$ and
     $y_\omega$. Thus we obtain a contradiction.

     \medskip

\textsc{ Step 2.} By iterating the result obtained in Step 1, one
can deduce that every pair of points in $\Gamma$ can be joined by
a discrete path of step at most $D$. Now it suffices to take the
finite set of all the elements in $\Gamma $ at distance at most
$D$ from $1$. By the previous statement, this is a set of
generators in $\Gamma$.

\me

\textbf{(2)}\quad According to (1) the group $\Gamma $ is finitely
generated; moreover it is easy to see that any word metric on
$\Gamma$ is bi-Lipschitz equivalent to $\dist$. Thus we may assume
that all asymptotic cones of some Cayley graph $\cgas$ are simply
connected.

The proof continues in the same spirit as for (1): there we dealt
with pairs of points and we had to ``fill the space between them''
with a discrete path composed of steps of bounded length. The set
of elements in the group with the above bound on their length gave
the finite set of generators. We can see a pair of points as an
image of the zero dimensional sphere $\sph^0$. If we go one
dimension up, instead of pairs of points we shall have loops in
the Cayley graph $\cgas$. These are nothing else than all the
relations in the group $\Gamma$ endowed with the finite generating
set $S$. To show that $\Gamma $ is finitely presented means to
show that an arbitrary loop in $\cgas$ can be ``filled'' with
loops of uniformly bounded length. ``Filled'' means here that by
putting a set of loops of uniformly bounded length one next to the
other one obtains a diagram having as boundary the initial
arbitrary loop. Then the (finite) set of words in the alphabet $S$
labelling loops in $\cgas$ of bounded length gives the set of
relations in the finite presentation.

More precisely, the argument comprises a two steps.

\medskip

\textsc{Step} 1. \quad We show that for every $M>0$ there exists
$N\in \N$ and $\ell_0$ such that every loop in $\cgas$ of length
$\ell \geq \ell_0$ can be filled by at most $N$ loops of length
$\leq \frac{\ell}{M}$.  This property is called the \textit{Loop
Division Property} in \cite{Dr1}; it is in fact equivalent to the
property that all asymptotic cones are simply connected (see
\cite{Dr1} for details).

The contrary of the Loop Division Property would imply that for
some $M>0$ there exists a sequence of loops $\cf_n $ of lengths
$\ell_n$ diverging to infinity such that for any set of at most
$n$ loops filling $\cf_n$, at least one of them has length larger
than $\frac{\ell_n}{M}$. We can see the loops $\cf_n$ as Lipschitz
maps $\cf_n : {\mathbb{S}}^1 \to \cgas$ with Lipschitz constant
$\frac{\ell_n}{2\pi}$.

In an asymptotic cone $\Gamma_\infty = \co{\Gamma ; x_n , \ell_n
}$ with $x_n$ on $\cf_n \left( {\mathbb{S}}^1 \right)$ the
sequence of loops $\cf_n$ defines a limit
$\frac{1}{2\pi}$--Lipschitz map $\cf : {\mathbb{S}}^1 \to
\Gamma_\infty$. Since $\Gamma_\infty$ is simply connected the map
$\cf$ can be extended to a continuous map $\bar \cf$ defined on
the unit disk $\mathbb{D}^2$. The uniform continuity of $\bar \cf$
implies that for a net on $\mathbb{D}^2$ of mesh $\delta$ small
enough its image by $\bar\cf$ is a union of $N$ ``squares'' of
perimeters at most $\frac{1}{2M}$ and filling $\cf$. Without loss
of generality it may be assumed that the edges of these
``squares'' are either geodesics or sub-arcs of $\cf$.

This implies that for $\omega$-almost every $n$ the loop $\cf_n$
can be filled by $N$ ``squares'' of perimeters at most
$\frac{\ell_n}{M}$, a contradiction.

     \medskip

\textsc{ Step 2.} By iterating the Loop Division Property obtained
in Step 1 we deduce that any loop in $\cgas$ can be filled by
loops of length at most $\ell_0$. There are finitely many loops of
length $\leq \ell_0$ up to left translations by elements in
$\Gamma$. The labels in the alphabet $S$ of these loops will be
the relations in the finite presentation of $\Gamma =<S>$.

\endproof

\begin{remark}
Simple connectedness of asymptotic cones implies much more than
the conclusion of Proposition \ref{sample}: it implies that the
group $\Gamma$ has polynomial Dehn function and linear filling
radius. See \cite{Dr1} and references therein.
\end{remark}

\subsection{Examples of asymptotic cones of groups}\label{exac}

All the groups considered below are finitely generated.

\me

\noindent  \textbf{(1)} A group is virtually nilpotent if and
    only if all its asymptotic cones are locally compact (\cite{Gr3}, \cite{Dr1}).

    \me

    \textit{Morality}: Outside the class of virtually nilpotent
    groups one should not expect the asymptotic cones to be
    locally compact.

    Moreover, in this case it was proved by P. Pansu in \cite{Pa} that
     all asymptotic cones are isometric to a graded Lie group canonically associated to $G$, as
    follows. Let $\mathrm{tor} (G)$ be the finite normal
    subgroup of $G$ generated by elements of finite order. The
    nilpotent group $\bar{G}=G/\mathrm{tor}(G)$ is without torsion, hence it can be embedded,
    according to \cite{Mal}, as a uniform lattice in a nilpotent Lie group. To this
    Lie group one canonically associates a graded Lie group, and it is this graded Lie group endowed with a Carnot-Caratheodory metric
     that is isometric to all
    asymptotic cones.

    If two virtually nilpotent groups are quasi-isometric, the
    graded Lie groups associated to them are not only bi-Lipschitz
    equivalent as usually for asymptotic cones, but moreover isomorphic. This points out new quasi-isometry invariants:
     the degree of nilpotency of $\bar{G}=G/\mathrm{tor}(G)$
    and the rank of each of the Abelian groups $\bar{G}^i/\bar{G}^{i+1}$, where $\bar{G}^i$ is the $i$-th group
     in the lower central series of $\bar{G}$.

    In particular, if a group $G$ is quasi-isometric to an Abelian group, then $G$ itself
    is virtually Abelian.

\me

\noindent  \textbf{(2)} A group is hyperbolic if and only if all
its asymptotic
    cones are real trees (\cite{Gr3}, \cite{Dr1}).

    Moreover, all asymptotic cones are isometric to a $2^{\aleph_0}$--universal
     real tree \cite{DP}.

\begin{remarks}
\begin{itemize}
    \item In the ``if'' part of the previous two statements as well as in
every similar statement in this paper, it is enough to take all
asymptotic cones for a \textit{fixed} ultrafilter.
    \item Note that Proposition \ref{sample}, (2), implies that hyperbolic
groups are finitely presented. This can also be obtained directly
from the definition, and in fact much more is known: the Dehn
function of hyperbolic groups is linear \cite{Gr2}.
\end{itemize}
\end{remarks}

\me

\noindent  \textbf{(3)} Let $G$ be a uniform lattice of isometries
of a symmetric space or Euclidean building of rank at least $2$.
Every asymptotic cone $\co{G;1,d}$ is a (non-discrete)
    Euclidean building \cite{KlL}.

    As for the question whether in this case all asymptotic cones
    are isometric or not, it turns out to be related to the
    Continuum Hypothesis (the hypothesis stating that there is no
    cardinal number between $\aleph_0$ and $2^{\aleph_0}$).

    Using a description of asymptotic cones in terms of fields and
    valuations (a similar description has been obtained independently by B. Leeb and A. Parreau
     \cite{Par}), Kramer, Shelah, Tent and Thomas have shown in \cite{KSTT} that:
     \begin{itemize}
        \item if the Continuum Hypothesis (CH) is not true then
        any uniform lattice in $SL(n,\R ),\, n\geq 3$, has
        $2^{2^{\aleph_0}}$ non-isometric asymptotic cones;
        \item if the CH is true then all asymptotic cones of a uniform lattice in $SL(n,\R ),\, n\geq
        3$, are isometric. Moreover, a finitely generated group
        has at most a continuum of non-isometric asymptotic cones.
     \end{itemize}

\me

\noindent  \textbf{(4)} In \cite{DS} can be found an example of
two-generated (and recursively presented - but not finitely
    presented) group with continuously many non-homeomorphic
    asymptotic cones. The construction is independent of CH.

\begin{ques} Can one characterize relatively hyperbolic groups
also in terms of asymptotic cones ?
\end{ques}

An answer to this question would give a better idea of how
relatively hyperbolic groups look like and also it might serve to
prove some rigidity result about relatively hyperbolic groups. For
instance, in the rigidity result of B. Kleiner and B.Leeb
\cite{KlL} the main ingredient is the description of asymptotic
cones of uniform lattices.

\section{Relatively hyperbolic groups: image from infinitely far away and rigidity}

\subsection{Tree-graded spaces and cut-points}

\begin{definition}[tree-graded spaces]\label{tree}
Let $\free$ be a complete geodesic metric space and let $\pp$ be a
collection of closed geodesic subsets of $\free$ (called
{\it{pieces}}) such that the following two properties are
satisfied:

\begin{enumerate}

\item[($T_1$)] Every two different pieces have at most one common
point.

\item[($T_2$)] Every simple geodesic triangle (a simple loop
composed of three geodesics) in $\free$ is contained in one piece.
\end{enumerate}

Then we say that the space $\free$ is {\em tree-graded with
respect to }$\pp$.
\end{definition}

\begin{figure}[!ht]
\centering
\includegraphics{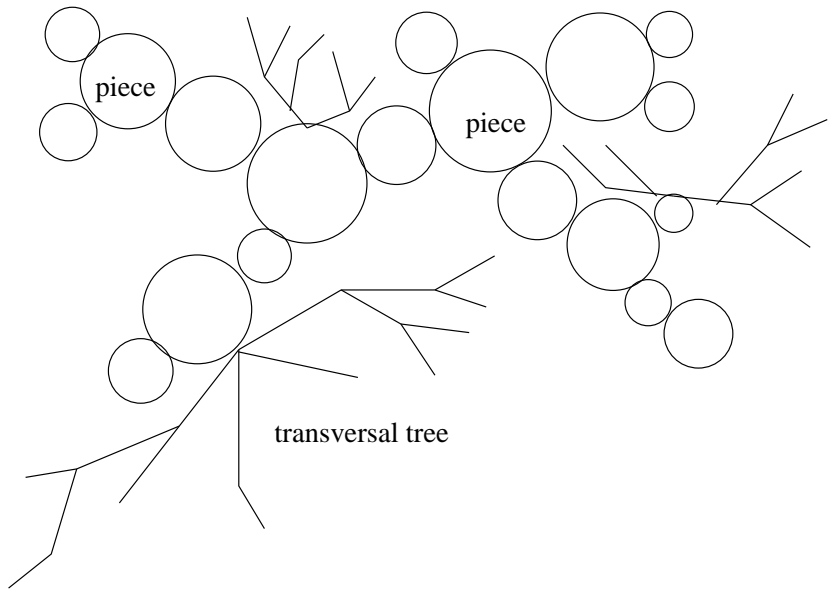}
\caption{A tree-graded space.} \label{fig6}
\end{figure}

Property ($T_2$) can be replaced by one of the following
properties:
\begin{enumerate}
\item[{\rm ($T_2'$)}]  For every topological arc $\cf:[0,d]\to \free$ and $t\in
[0,d]$, let $\cf[t-a,t+b]$ be a maximal sub-arc of $\cf$
containing $\cf (t)$ and contained in one piece. Then every other
topological arc with the same endpoints as $\cf$ must contain the
points $\cf (t-a)$ and $\cf (t+b)$.
\begin{figure}[!ht]
\centering
\unitlength .85mm 
\linethickness{0.4pt}
\ifx\plotpoint\undefined\newsavebox{\plotpoint}\fi 
\begin{picture}(128.04,45.33)(0,0)
\put(93.47,23.86){\line(0,1){.97}}
\put(93.44,24.83){\line(0,1){.968}}
\put(93.38,25.8){\line(0,1){.964}}
\multiput(93.27,26.77)(-.0306,.19171){5}{\line(0,1){.19171}}
\multiput(93.11,27.73)(-.0327,.15844){6}{\line(0,1){.15844}}
\multiput(92.92,28.68)(-.02987,.1176){8}{\line(0,1){.1176}}
\multiput(92.68,29.62)(-.03125,.10323){9}{\line(0,1){.10323}}
\multiput(92.4,30.55)(-.0323,.09154){10}{\line(0,1){.09154}}
\multiput(92.08,31.46)(-.03309,.0818){11}{\line(0,1){.0818}}
\multiput(91.71,32.36)(-.03369,.07354){12}{\line(0,1){.07354}}
\multiput(91.31,33.24)(-.031701,.061663){14}{\line(0,1){.061663}}
\multiput(90.86,34.11)(-.032159,.056156){15}{\line(0,1){.056156}}
\multiput(90.38,34.95)(-.032499,.051229){16}{\line(0,1){.051229}}
\multiput(89.86,35.77)(-.032736,.046783){17}{\line(0,1){.046783}}
\multiput(89.3,36.56)(-.032883,.042741){18}{\line(0,1){.042741}}
\multiput(88.71,37.33)(-.032951,.039042){19}{\line(0,1){.039042}}
\multiput(88.09,38.08)(-.032948,.035637){20}{\line(0,1){.035637}}
\multiput(87.43,38.79)(-.032882,.032486){21}{\line(-1,0){.032882}}
\multiput(86.74,39.47)(-.036033,.032515){20}{\line(-1,0){.036033}}
\multiput(86.02,40.12)(-.039438,.032476){19}{\line(-1,0){.039438}}
\multiput(85.27,40.74)(-.043136,.032364){18}{\line(-1,0){.043136}}
\multiput(84.49,41.32)(-.047176,.032167){17}{\line(-1,0){.047176}}
\multiput(83.69,41.87)(-.051618,.031877){16}{\line(-1,0){.051618}}
\multiput(82.86,42.38)(-.060579,.033726){14}{\line(-1,0){.060579}}
\multiput(82.01,42.85)(-.066814,.033334){13}{\line(-1,0){.066814}}
\multiput(81.15,43.28)(-.07394,.0328){12}{\line(-1,0){.07394}}
\multiput(80.26,43.68)(-.0822,.0321){11}{\line(-1,0){.0822}}
\multiput(79.35,44.03)(-.09192,.03119){10}{\line(-1,0){.09192}}
\multiput(78.44,44.34)(-.1036,.03){9}{\line(-1,0){.1036}}
\multiput(77.5,44.61)(-.1348,.03251){7}{\line(-1,0){.1348}}
\multiput(76.56,44.84)(-.15883,.03078){6}{\line(-1,0){.15883}}
\multiput(75.61,45.02)(-.19207,.02828){5}{\line(-1,0){.19207}}
\put(74.65,45.17){\line(-1,0){.966}}
\put(73.68,45.26){\line(-1,0){.969}}
\put(72.71,45.32){\line(-1,0){.971}}
\put(71.74,45.33){\line(-1,0){.97}}
\put(70.77,45.29){\line(-1,0){.968}}
\multiput(69.8,45.22)(-.2408,-.0303){4}{\line(-1,0){.2408}}
\multiput(68.84,45.09)(-.19133,-.03292){5}{\line(-1,0){.19133}}
\multiput(67.88,44.93)(-.13546,-.02967){7}{\line(-1,0){.13546}}
\multiput(66.93,44.72)(-.11723,-.03129){8}{\line(-1,0){.11723}}
\multiput(66,44.47)(-.10284,-.0325){9}{\line(-1,0){.10284}}
\multiput(65.07,44.18)(-.09114,-.0334){10}{\line(-1,0){.09114}}
\multiput(64.16,43.85)(-.07461,-.03124){12}{\line(-1,0){.07461}}
\multiput(63.26,43.47)(-.067501,-.031922){13}{\line(-1,0){.067501}}
\multiput(62.39,43.06)(-.061275,-.032445){14}{\line(-1,0){.061275}}
\multiput(61.53,42.6)(-.055762,-.032836){15}{\line(-1,0){.055762}}
\multiput(60.69,42.11)(-.050832,-.033116){16}{\line(-1,0){.050832}}
\multiput(59.88,41.58)(-.046384,-.033299){17}{\line(-1,0){.046384}}
\multiput(59.09,41.01)(-.04234,-.033398){18}{\line(-1,0){.04234}}
\multiput(58.33,40.41)(-.03864,-.033421){19}{\line(-1,0){.03864}}
\multiput(57.59,39.78)(-.035235,-.033377){20}{\line(-1,0){.035235}}
\multiput(56.89,39.11)(-.03369,-.034936){20}{\line(0,-1){.034936}}
\multiput(56.22,38.41)(-.032077,-.036423){20}{\line(0,-1){.036423}}
\multiput(55.57,37.68)(-.031997,-.039828){19}{\line(0,-1){.039828}}
\multiput(54.97,36.93)(-.033712,-.046085){17}{\line(0,-1){.046085}}
\multiput(54.39,36.14)(-.033569,-.050534){16}{\line(0,-1){.050534}}
\multiput(53.86,35.33)(-.033333,-.055467){15}{\line(0,-1){.055467}}
\multiput(53.36,34.5)(-.032991,-.060983){14}{\line(0,-1){.060983}}
\multiput(52.89,33.65)(-.032523,-.067213){13}{\line(0,-1){.067213}}
\multiput(52.47,32.77)(-.03191,-.07433){12}{\line(0,-1){.07433}}
\multiput(52.09,31.88)(-.0311,-.08258){11}{\line(0,-1){.08258}}
\multiput(51.75,30.97)(-.03341,-.10255){9}{\line(0,-1){.10255}}
\multiput(51.45,30.05)(-.03234,-.11695){8}{\line(0,-1){.11695}}
\multiput(51.19,29.12)(-.03088,-.13519){7}{\line(0,-1){.13519}}
\multiput(50.97,28.17)(-.02886,-.15919){6}{\line(0,-1){.15919}}
\multiput(50.8,27.21)(-.0324,-.2405){4}{\line(0,-1){.2405}}
\put(50.67,26.25){\line(0,-1){.967}}
\put(50.58,25.28){\line(0,-1){3.877}}
\multiput(50.68,21.41)(.0332,-.2404){4}{\line(0,-1){.2404}}
\multiput(50.81,20.45)(.02936,-.15909){6}{\line(0,-1){.15909}}
\multiput(50.98,19.49)(.0313,-.13509){7}{\line(0,-1){.13509}}
\multiput(51.2,18.55)(.03271,-.11684){8}{\line(0,-1){.11684}}
\multiput(51.47,17.61)(.03374,-.10244){9}{\line(0,-1){.10244}}
\multiput(51.77,16.69)(.03137,-.08248){11}{\line(0,-1){.08248}}
\multiput(52.11,15.78)(.03214,-.07423){12}{\line(0,-1){.07423}}
\multiput(52.5,14.89)(.032736,-.067109){13}{\line(0,-1){.067109}}
\multiput(52.93,14.02)(.033184,-.060878){14}{\line(0,-1){.060878}}
\multiput(53.39,13.17)(.033509,-.055361){15}{\line(0,-1){.055361}}
\multiput(53.89,12.34)(.033729,-.050428){16}{\line(0,-1){.050428}}
\multiput(54.43,11.53)(.031977,-.043423){18}{\line(0,-1){.043423}}
\multiput(55.01,10.75)(.032123,-.039726){19}{\line(0,-1){.039726}}
\multiput(55.62,9.99)(.032192,-.036322){20}{\line(0,-1){.036322}}
\multiput(56.26,9.27)(.032191,-.033171){21}{\line(0,-1){.033171}}
\multiput(56.94,8.57)(.035341,-.033265){20}{\line(1,0){.035341}}
\multiput(57.65,7.91)(.038746,-.033298){19}{\line(1,0){.038746}}
\multiput(58.38,7.27)(.042446,-.033263){18}{\line(1,0){.042446}}
\multiput(59.15,6.67)(.046489,-.033152){17}{\line(1,0){.046489}}
\multiput(59.94,6.11)(.050937,-.032955){16}{\line(1,0){.050937}}
\multiput(60.75,5.58)(.055866,-.03266){15}{\line(1,0){.055866}}
\multiput(61.59,5.09)(.061377,-.032251){14}{\line(1,0){.061377}}
\multiput(62.45,4.64)(.067601,-.031708){13}{\line(1,0){.067601}}
\multiput(63.33,4.23)(.07471,-.031){12}{\line(1,0){.07471}}
\multiput(64.22,3.86)(.09125,-.03311){10}{\line(1,0){.09125}}
\multiput(65.14,3.53)(.10294,-.03217){9}{\line(1,0){.10294}}
\multiput(66.06,3.24)(.11733,-.03092){8}{\line(1,0){.11733}}
\multiput(67,2.99)(.13555,-.02924){7}{\line(1,0){.13555}}
\multiput(67.95,2.78)(.19143,-.03231){5}{\line(1,0){.19143}}
\put(68.91,2.62){\line(1,0){.963}}
\put(69.87,2.5){\line(1,0){.968}}
\put(70.84,2.43){\line(1,0){.97}}
\put(71.81,2.4){\line(1,0){.971}}
\put(72.78,2.41){\line(1,0){.969}}
\put(73.75,2.47){\line(1,0){.965}}
\multiput(74.71,2.57)(.19198,.02889){5}{\line(1,0){.19198}}
\multiput(75.67,2.72)(.15873,.03128){6}{\line(1,0){.15873}}
\multiput(76.63,2.9)(.1347,.03294){7}{\line(1,0){.1347}}
\multiput(77.57,3.13)(.1035,.03033){9}{\line(1,0){.1035}}
\multiput(78.5,3.41)(.09182,.03148){10}{\line(1,0){.09182}}
\multiput(79.42,3.72)(.0821,.03236){11}{\line(1,0){.0821}}
\multiput(80.32,4.08)(.07384,.03304){12}{\line(1,0){.07384}}
\multiput(81.21,4.47)(.066708,.033545){13}{\line(1,0){.066708}}
\multiput(82.07,4.91)(.056441,.031657){15}{\line(1,0){.056441}}
\multiput(82.92,5.38)(.051517,.03204){16}{\line(1,0){.051517}}
\multiput(83.75,5.9)(.047074,.032317){17}{\line(1,0){.047074}}
\multiput(84.55,6.45)(.043033,.0325){18}{\line(1,0){.043033}}
\multiput(85.32,7.03)(.039335,.032601){19}{\line(1,0){.039335}}
\multiput(86.07,7.65)(.03593,.032629){20}{\line(1,0){.03593}}
\multiput(86.79,8.3)(.032779,.03259){21}{\line(1,0){.032779}}
\multiput(87.47,8.99)(.032835,.035741){20}{\line(0,1){.035741}}
\multiput(88.13,9.7)(.032827,.039146){19}{\line(0,1){.039146}}
\multiput(88.75,10.45)(.032748,.042845){18}{\line(0,1){.042845}}
\multiput(89.34,11.22)(.032587,.046887){17}{\line(0,1){.046887}}
\multiput(89.9,12.02)(.032336,.051332){16}{\line(0,1){.051332}}
\multiput(90.42,12.84)(.031982,.056257){15}{\line(0,1){.056257}}
\multiput(90.9,13.68)(.031506,.061763){14}{\line(0,1){.061763}}
\multiput(91.34,14.54)(.03346,.07364){12}{\line(0,1){.07364}}
\multiput(91.74,15.43)(.03283,.08191){11}{\line(0,1){.08191}}
\multiput(92.1,16.33)(.03201,.09164){10}{\line(0,1){.09164}}
\multiput(92.42,17.25)(.03092,.10333){9}{\line(0,1){.10333}}
\multiput(92.7,18.18)(.03371,.13451){7}{\line(0,1){.13451}}
\multiput(92.93,19.12)(.0322,.15854){6}{\line(0,1){.15854}}
\multiput(93.13,20.07)(.02999,.19181){5}{\line(0,1){.19181}}
\put(93.28,21.03){\line(0,1){.965}}
\put(93.38,21.99){\line(0,1){1.871}}
\put(57,21.86){\circle*{1.8}} \put(84.75,22.11){\circle*{.71}}
\put(85,22.36){\circle*{1}}
\qbezier(14.5,19.11)(33.88,15.11)(56.75,22.11)
\qbezier(57,22.86)(70.75,28.36)(84.5,22.86)
\qbezier(84.75,22.86)(110.88,15.74)(127.5,23.11)
\put(14.5,19.11){\circle*{1.58}}
\put(127.25,23.11){\circle*{1.58}} \thicklines
\qbezier(14.25,19.36)(25.75,22.86)(35.25,18.36)
\qbezier(57,21.86)(72.88,43.24)(85.25,22.11)
\put(22.75,15.36){\makebox(0,0)[cc]{$\cf$}}
\put(50.61,22.76){\makebox(0,0)[cc]{$\cf(t-a)$}}
\put(79.81,20.55){\makebox(0,0)[cc]{$\cf(t+b)$}}
\put(84.93,22.6){\circle*{1.69}}
\qbezier(35.32,18.08)(47.35,11.98)(56.66,21.86)
\qbezier(85.04,22.7)(103.43,-5.99)(127.29,22.91)
\end{picture}
\caption{Property ($T_2'$).} \label{fig7}
\end{figure}
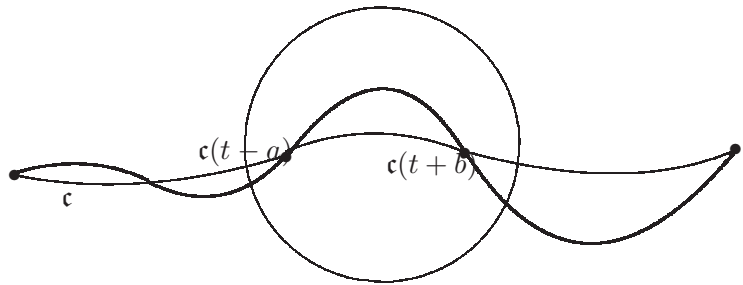
\item[{\rm ($T_2''$)}] Every
simple loop in $\free$ is contained in one piece.
\end{enumerate}

\begin{remark}[\cite{DS}]\label{geodsub}
If one replaces property $(T_2)$ by the stronger property
$(T_2'')$ in the definition of a tree-graded space then one can
weaken the condition on $\calp$ and ask only that each set in
$\calp$ is path-connected.
\end{remark}

The structure of tree-graded space appears naturally as soon as a
space has a cut-point, as shown by the following result.

\begin{proposition}[\cite{DS}, Section $\S 2.4$]\label{cutp}
Let $X$ be a complete geodesic metric space containing at least
two points and let $\calc$ be a non-empty set of global cut-points
in $X$.
\begin{itemize}
\item[(a)] There exists a largest (in an appropriate sense)
collection $\calp$ of subsets of $X$ such that $X$ is tree-graded
with respect to $\calp$ and such that any piece in $\calp$ is
either a singleton or a set with no global cut-point in $\calc$.

Moreover the intersection of any two distinct pieces from $\calp$
is either empty or a point in~$\calc$.
\item[(b)] If $\calc=X$ then all pieces in $\calp$ are either
singletons or sets without cut-point. In particular this is true
if $X$ is a homogeneous space with a cut-point.
\end{itemize}
\end{proposition}

We should point out here that a systematic approach to spaces
having cut-points had already been done by B. Bowditch when
working at the Bestvina-Mess conjecture (stating that the boundary
at infinity of a one-ended hyperbolic group is locally connected,
or equivalently that it has no global cut-point). Thus, Bowditch
considered a topological space with a cut-point and such that the
group of its homeomorphisms acts on it with dense orbits (having
in mind the boundary at infinity of a group), and he showed that
such a space projects on a tree. For details and references see
\cite{Bow1}.

\bi

\noindent \textit{Properties of tree-graded spaces}:
\begin{enumerate}
    \item If all the pieces are real trees then $\free$ is a real
    tree.
    \item For every $x\in \free $ we define the set
    $\mathcal{T}_x$ to be the set of points $y\in \free$ which can be joined to $x$ by a topological arc
    intersecting each piece in at most a point. For every $x$ the
    set $\mathcal{T}_x$ is a real tree and a closed subset in
    $\free$. For every $y\in \mathcal{T}_x$,
    $\mathcal{T}_y=\mathcal{T}_x$. We call such a tree \textit{transversal
    tree}.
    \item Every point of intersection of two distinct pieces as well as
    every point in a non-trivial transversal tree is a cut-point for $\free$.
    \item Every path-connected subset \textit{without cut-points}
    is contained in a piece.
    \item For every point $x$ outside a piece $M$ there exists a unique
    point on $M$ minimizing the distance to $x$. This allows to
    define a projection map from $\free $ to $M$.
    \item Every path-connected subset intersecting a piece $M$ in at
    most one point projects onto the piece $M$ in a unique point.
    \item Suppose that there exists $\epsilon >0$ such that every
    loop of length at most $\epsilon$ and contained in a piece is
    contractible. Then $\pi_1 (\free )$ coincides with the free product of the fundamental groups
     of the pieces, $\displaystyle \ast_{M\in \pp }\pi_1 (M)$.

     We note that the hypothesis that short loops are contractible
     in any piece is a necessary condition, as shown by the
     example of the Hawaiian earring.
    \item Let $(\free , \pp )$ be a tree-graded space. If $\phi$
    is a homeomorphism from $\free $ to another geodesic metric space
    $X$, then $X$ is tree-graded with respect to the collection of
    pieces $\{\phi(M)\mid M\in \pp \}$.
    \item Let $(\free ,\pp )$ and $(\free' , \pp')$ be two
    tree-graded spaces with all pieces without cut-points. Every
    homeomorphism $\phi :\free \to \free'$ sends each piece onto a
    piece.
\end{enumerate}

For proofs of these properties and other properties of tree-graded
spaces see \cite{DS}.

\me

\subsection{The characterization of relatively hyperbolic groups in
terms of asymptotic cones}

\begin{theorem}[Dru\c{t}u-Osin-Sapir\cite{DS}]\label{thrhc}
A finitely generated group $G$ is hyperbolic relatively to a
finite family $\{ H_1,...,H_n \}$ of finitely generated subgroups
if and only if every asymptotic cone $\Con_\omega(G;1,d)$ is
tree-graded with respect to the collection of pieces
$$
\pp = \left\{ \lio{g_n H_i}\mid (g_n)\mbox{ sequence in G}, i\in
\{ 1,2,\dots ,n\} \right\}\, .
$$
\end{theorem}

The ``only if'' part is proven by D. Osin and M. Sapir in the
Appendix of \cite{DS}. The ``if'' part is proven in \cite{DS}.
Note that for the ``if'' part one does not need to ask that the
peripheral subgroups are finitely generated. It follows
immediately from the fact that the limit sets of their left
cosets, which are isometric to their asymptotic cones with the
metric induced from $G$, are geodesic, since they are pieces in a
tree-graded space. It remains to apply Proposition \ref{sample},
(1).

Also, Proposition \ref{sample}, (2), and property 7 of tree-graded
spaces implies that if $H_i$ all have simply connected asymptotic
cones then $G$ is finitely presented. On the other hand, from the
equivalent definition of relative hyperbolicity given by D. Osin
in \cite{Osin} it follows that the same is true if all $H_i$ are
finitely presented (which is a weaker hypothesis than the
previous).

In particular Theorem \ref{thrhc} is true for $G=\Gamma $ a
non-uniform lattice in rank one and $\{ \Gamma_{\alpha_1},\dots ,
\Gamma_{\alpha_m}\}$ its cusp subgroups. Thus the image of the
space $X_0$ in Figure \ref{fig4} seen from infinitely far away is
a homogeneous version of the Figure \ref{fig6}. It is not
difficult to show that for this particular tree-graded space each
transversal tree is in fact a $2^{\aleph_0}$--universal real tree.

A straightforward consequence of Theorem \ref{thrhc} is the
following.

\begin{cor}\label{cor} If a group $G$ is hyperbolic relative to
$\{H_1,\dots ,H_m\}$ and if each $H_i$ is hyperbolic relative to a
collection of subgroups $\{ H_i^{1}, \dots , H_i^{n_i} \}$ then
$G$ is hyperbolic relative to $\{ H_i^{j} \mid i\in \{1,\dots ,m
\},\, j\in \{1,\dots ,n_i \}\}$. \end{cor}

\begin{remark}\label{coment}
This process may not terminate: for instance if $G$ is a free
group and $H= \langle h \rangle$ is a cyclic subgroup, one can
consider $H_n=\langle h^{2^n} \rangle$, $G$ is hyperbolic relative
to $\{ H_n\}$ and $H_n$ is hyperbolic relative to $\{ H_{n+1}\}$.
Still, in this situation there exists a terminal point: $G$
hyperbolic relative to $\{ 1\}$.

In general, a terminal point would be a family $\{H_1,\dots
,H_m\}$ of peripheral subgroups relative to which the ambient
group $G$ is hyperbolic and such that no $H_i$ is relatively
hyperbolic. Such a family may not exist for an arbitrary
relatively hyperbolic group. Indeed, the example of inaccessible
group constructed by Dunwoody in \cite{Du2} is also an example of
relatively hyperbolic group such that every list of peripheral
subgroups must contain a relatively hyperbolic subgroup (the
argument showing this can be found in \cite{BDM}). See Question
\ref{last} in Section \ref{soq}.
\end{remark}

\bigskip

Theorem \ref{thrhc} and properties (2) and (9) of tree-graded
spaces suggest that the ``good objects'' for a rigidity theory for
relatively hyperbolic groups are the finitely generated groups
such that all their asymptotic cones are without cut-points. We
call such groups \textit{asymptotically without cut-points}. To
avoid trivial cases and different technical complications we also
assume that finite groups are not asymptotically without
cut-points.

\begin{remark}\label{roneend}
A group asymptotically without cut-points is one-ended. This
follows from Stallings' Ends Theorem stating that a finitely
generated group splits as a free product or HNN-extension with
finite amalgamation if and only if it has more than one end
\cite{Sta}. The converse is not true: the asymptotic cones of any
hyperbolic group are $\R$--trees, and there are one-ended
hyperbolic groups (uniform lattices in $\hip^3_\R$ for instance).
\end{remark}

\me

\noindent \textit{Examples of groups asymptotically without
cut-points}:
\begin{enumerate}
    \item products $G=G_1\times G_2$, where $G_1$ and $G_2$ are
    infinite groups; this follows from the fact that any
    asymptotic cone of $G$ is a product of asymptotic cones
    of $G_1$ and of $G_2$, and the latter are geodesic spaces;
    \item uniform lattices in symmetric spaces/Euclidean buildings
    of rank at least two; this is because their asymptotic cones
    are non-discrete Euclidean buildings of rank at least two, and
    these do not have cut-points \cite{KlL};
    \item groups with elements of infinite order in the center, not virtually cyclic
    (\cite{DS}, see also paragraph 6.1 in the present paper);
    \item groups satisfying an identity (a law), not virtually cyclic
    (\cite{DS}, see also paragraph 6.2 in the present paper).

    We recall what \textit{satisfying an identity (a law)} means for a group.
    Let $w(x_1,\dots ,x_n)$ be a non-trivial reduced word in the $n$ letters $x_1,\dots
    ,x_n$ and their inverses. \textit{Reduced} means that all sequences of type $x\, x^{-1}$ are deleted.
     The group $G$ satisfies the identity
    $w(x_1,\dots ,x_n)=1$ if the equality is satisfied in $G$
    whenever replacing $x_1,\dots ,x_n$ with arbitrary elements in
    $G$.

    Examples of such groups:
    \begin{itemize}
        \item Abelian groups: here $w=x_1 x_2 x_1^{-1} x_2^{-1}$;
        \item more generally solvable groups of class at most $m\in
        \N$;
        \item free Burnside groups. We recall that the free Burnside group
        $B(n,m)$ is the group with $n$ generators satisfying the
        identity $x^m=1$ and all the relations that can be
        obtained from this identity (and no other). A rigorous way to define it is to say that it is the quotient
         of $\free_n$ by its normal subgroup generated by all elements of the form $f^m,\, f\in \free_m$.
          It is known that these groups are infinite for $m$ large enough
          (see \cite{Ad}, \cite{Olsh:book}, \cite{I}, \cite{Ly}, \cite{DG} and references
        therein).
        \item uniformly amenable groups, not virtually cyclic.

A discrete group $G$ is \textit{uniformly amenable} if there
exists a function ${\mathfrak{C}} : (0,1) \times \N \to \N$ such
that for every finite subset $K$ of $G$ and every $\epsilon \in
 (0,1)$ there exists a finite subset $F\subset G$ satisfying:
 \begin{itemize}
    \item[(i)] card$\, F\leq {\mathfrak{C}} (\epsilon ,
    \mathrm{card}\, K)$;
    \item[(ii)] card$\, KF < (1+\epsilon ) \mathrm{card}\, F$.
 \end{itemize}

For details on this notion see \cite{Kel}, \cite{Boz} and
\cite{Wys}. In \cite{DS} it is shown that a uniformly amenable
group always satisfies a law.

\end{itemize}
\end{enumerate}

\subsection{Rigidity of relatively hyperbolic
groups}\label{sectrig}

\begin{theorem}[\cite{DS}, \cite{BDM}] \label{s}
Let $G$ be a finitely generated group that is hyperbolic relative
to its subgroups $H_1, ...,H_m$, and let $S$ be a finitely
generated group that is not relatively hyperbolic with respect to
any finite collection of proper subgroups.

Then the image of $S$ under any $(L,C)$--quasi-isometric embedding
$S\to G$ is in the $M$--tubular neighborhood of a coset $gH_i$,
$g\in G, i=1,...,m$, where $M$ depends only on $L, C, G$ and $H_1,
...,H_m$.
\end{theorem}

\begin{remark}
In \cite[$\S 3$]{PW} Theorem \ref{s} is proven for $G$
     a fundamental group of a graph of groups with finite edge
groups and $S$ a one-ended group.

One cannot hope however to weaken the hypothesis of Theorem
\ref{s} to ``$S$ a one ended group''. For instance non-uniform
lattices in $\hip^3_\R$ are one-ended groups on one hand and
hyperbolic relative to their cusp subgroups on the other. Thus,
they are quasi-isometrically embedded into themselves and not
uniformly near a left coset of a cusp subgroup.
\end{remark}

\begin{cor}\label{ssgp}
Let $G$ be a finitely generated group that is hyperbolic relative
to its subgroups $H_1, ...,H_m$, and let $S$ be an undistorted
subgroup of $G$ that is not relatively hyperbolic with respect to
any finite collection of proper subgroups. Then $S$ is contained
in $H_i^g$ for some $g\in G$ and $ i\in \{ 1,...,m\} $.
\end{cor}

As the proofs of the Theorems of R. Schwartz show, a rigidity
result such as Theorem \ref{s} can be used to get a result on q.i.
completeness. Indeed, this can also be done in this case.

\begin{theorem}[\cite{DS}, , \cite{BDM}] \label{rigrh} Let
$G$ be a finitely generated group hyperbolic relative to $\{ H_1,
...,H_m\}$. Suppose that all the subgroups $H_i$, $i=1,...,m$, are
 not relatively hyperbolic with respect to
any finite collection of proper subgroups.

Let $\Lambda $ be a finitely generated group that is
quasi-isometric to $G$. Then $\Lambda $ is hyperbolic relative to
a finite collection of subgroups $S_1,...,S_n$ each of which is
quasi-isometric to one of the subgroups $H_1,...,H_m$.
\end{theorem}

\begin{remarks}
\begin{itemize}
    \item[(a)] The number of subgroups relative to which the group
    is hyperbolic is \textit{not} a quasi-isometry invariant. This
    can be seen in the example of a finite covering $M \to N$ of a
    finite volume non-compact hyperbolic $3$-manifold by another.
    The group $\Gamma_M=\pi_1 (M)$ is a finite index subgroup of
    $\Gamma_N=\pi_1 (N)$, so they are quasi-isometric. On the
    other hand, the number of cusp subgroups of $\Gamma_M$ can be
    larger than the number of cusp subgroups of~$\Gamma_N$.
    \item[(b)] Particular cases of Theorem \ref{rigrh} follow from the
    results in \cite{Sch}, \cite{KaL1}, \cite{KaL2}, \cite{PW}.
\end{itemize}
\end{remarks}

In view of Theorems \ref{s} and \ref{rigrh} it becomes interesting
to provide a list of

\me

\noindent \textit{Examples of groups that are not relatively
hyperbolic:}
\begin{enumerate}
    \item groups having one asymptotic cone without cut-point;
    that such groups are not relatively hyperbolic follows from Theorem \ref{thrhc};
    \item groups without free
    non-Abelian subgroups and not virtually cyclic; this follows from the fact that relatively
hyperbolic groups that are not virtually cyclic contain free
non-Abelian subgroups.

This class of groups contains the amenable groups not virtually
cyclic, but it is strictly larger than that class. See Remark
\ref{rvn} for more details.
    \item groups with infinite
    center and not virtually cyclic; indeed, if a group is  not virtually cyclic and it is
     hyperbolic relative to proper subgroups then its center is finite.
\end{enumerate}

Other examples can be found in Section \ref{scutnrh}.

\bigskip

\noindent \textit{Outline of proof of Theorem \ref{rigrh}}.

\me

\noindent \textbf{I.} Let $\q :\Lambda \to G$ and $\bar{\q} :G \to
\Lambda $ be two $(L,C)$--quasi-isometries quasi-converse to each
other. Lemma \ref{quasia} implies that using them one can
construct a quasi-action quasi-transitive and of finite kernel of
$\Lambda$ on the Cayley graph of $G$. For simplicity we denote by
$\aaa$ the set $\{ gH_i \mid g\in G, i=1,2,\dots ,m \}$. By
quasi-transitivity, for every left coset $A\in \aaa$ and every
point $g$ in it, there exists $\lambda \in \Lambda$ such that
$\q_\lambda (g) \in B(1,C_1)$, where $C_1=C_1(L,C)$. On the other
hand, by Theorem \ref{s}, $\q_\lambda (A)$ is contained in the
$M$--tubular neighborhood of another left coset $A'\in \aaa$,
where $M=M(L,C,G)$. It follows that $A'$ intersects $B(1,C_1+M)$.

We conclude that the finite set $\{A_1 ... A_k\}$ of left cosets
intersecting $B(1,C_1+M)$, is in some sense a set of
representatives for $A\in \aaa$, therefore we will diminish it (in
the next step), and keep only one orbit representative by
conjugacy class.

\me

\noindent \textbf{II.} First we note that if for some $\lambda$ in
$ \Lambda $, $A$ and $B$ in $\aaa$, and $M>0$ we have that
$\q_\lambda (A) \subset \nn_M (B)$, then $\dist_H (\q_\lambda (A)
, B)\leq M'$ for some $M'=M'(L,C)$. This follows by applying the
same result to $\q_{\lambda^{-1}}(B)$ and from the fact that two
left cosets cannot be at finite Hausdorff distance one from the
other unless they coincide.

Now we consider the equivalence relation in $\aaa$
$$
A \sim B \Leftrightarrow \exists \lambda \in \Lambda,\, \exists
M>0 \mbox{ such that } \q_\lambda (A)\subset \nn_M (B)\, .
$$

In the set $\{ A_1,\dots , A_k\}$ we select one representative in
each equivalence class and obtain thus a possibly smaller set $\{
B_1,\dots , B_n\}$. Also, for every $A_i \sim A_j$ we fix
$\lambda_{ij}$ such that $\q_{\lambda_{ij}} (A_i)\subset \nn_M
(A_j)$, and we consider $K_0=\max_{i,j} \dist
(\q_{\lambda_{ij}}(1),1)$.

 \me

\noindent \textbf{III.} We define for each $A\in \aaa$ the
subgroup in $\Lambda$
$$
\mathrm{Stab}_{M'}(A)=\{ \lambda\in \Lambda \: ;\:
\dist_H(\q_\lambda (A) , A)\leq M'\}\, .
$$

Using the arguments in \textbf{I} and the choice made in
\textbf{II} it is not difficult to show that for every $A\in \aaa
$, $\mathrm{Stab}_{M'}(A)$ acts $C_2$--quasi-transitively on $A$,
in the sense that every orbit of a point under the action of the
group contains $A$ in its $C_2$--tubular neighborhood. Here $C_2$
is a constant which is computed by means of $K_0$.

This in particular implies that $\dist_H (\bar{\q}(B_i),
\mathrm{Stab}_{M'}(B_i))\leq \kappa $ for some constant $\kappa =
\kappa (L,C,C_2)$. The last statement together with the argument
in \textbf{I} imply that for every $A\in \aaa$ there exists
$\lambda \in \Lambda $ and $i\in \{ 1,2,\dots ,n\}$ such that
$$
\dist_H (\bar{\q}(A), \lambda \mathrm{Stab}_{M'}(B_i) )\leq \chi\,
,
$$ for some $\chi =\chi (L,C,\kappa )$.

\me

\noindent \textbf{IV.} Now we have the following sequence of
implications. $G$ is hyperbolic relative to $\{ H_1,\dots , H_m
\}$ $\Rightarrow$ (by Theorem \ref{thrhc}) every asymptotic cone
$\co{G;1,d}$ is tree-graded with set of pieces $\{ \lio{A_n}\}\:
;\: A_n \in \aaa \}$ $\Rightarrow $ (by Property 8 of tree-graded
spaces and Property 4 of asymptotic cones) every asymptotic cone
$\co{\Lambda ;1,d}$ is tree-graded with set of pieces $\{
\lio{\bar{\q}(A_n)}\}\: ;\: A_n \in \aaa \}$ $\Rightarrow$ (by the
last statement in \textbf{III}) every asymptotic cone $\co{\Lambda
;1,d}$ is tree-graded with set of pieces $\{ \lio{\lambda_n
\mathrm{Stab}_{M'}(B_i)}\: ;\: \lambda_n \in \Lambda ,\, i\in \{
1,2,\dots ,n\} \}$ $\Rightarrow $ (again by Theorem \ref{thrhc})
$\Lambda $ is hyperbolic relative to $\{
\mathrm{Stab}_{M'}(B_1),\dots , \mathrm{Stab}_{M'}(B_n) \}$.
\hspace*{\fill}$\square$

\bigskip

Finally, it turns out that the whole class of groups hyperbolic
relative to proper subgroups is q.i. complete.

\begin{theorem}[relative hyperbolicity is q.i. invariant \cite{Dr3}]\label{iqir}
Let $G$ be a group hyperbolic relative to a family of subgroups
$H_1,...,H_n$. If a group $G'$ is quasi-isometric to $G$ then $G'$
is hyperbolic relative to $H_1',...,H_m'$, where each $H_i'$ can
be embedded quasi-isometrically in $H_j$ for some $j=j(i)\in
\{1,2,...,n\}$.
\end{theorem}

Note that Theorem \ref{rigrh} does not imply Theorem \ref{iqir}
because there exist relatively hyperbolic groups that have no list
of peripheral subgroups composed uniquely of subgroups not
relatively hyperbolic (see Remark \ref{coment}).

Note also that in the full generality assumed in Theorem
\ref{iqir}, the stronger conclusion that each subgroup $H_i'$ is
quasi-isometric to some subgroup $H_j$ cannot hold. This can be
seen for instance when $G=G'=A*B*C$, with $G$ hyperbolic relative
to $\{A*B,C\}$ and $G'$ hyperbolic relative to $\{A,B,C\}$.

The proof of Theorem \ref{iqir} has an outline completely
different from the one of Theorem \ref{rigrh}. Its main ingredient
is not (and cannot be) a quasi-isometric embedding rigidity result
as Theorem \ref{s}. But it relies on some new geometric ways to
define relative hyperbolicity.

\subsection{More rigidity of relatively hyperbolic groups: outer
automorphisms}\label{sectout}

The group of outer automorphisms of a group $G$ is the quotient
group $Out(G)=Aut(G)/Inn(G)$, where $Inn(G)$ is the normal
subgroup of automorphisms $c_g$ given by the conjugacy with an
element $g\in G$. The group $Inn(G)$ is called \textit{the group
of inner automorphisms}.

We recall that in the case of hyperbolic groups the following
result is known.

\begin{theorem}[\cite{Pau}]\label{out}
\begin{itemize}
    \item[(1)] Let $G $ be a hyperbolic group. If $Out(G )$ is infinite then $G$
    acts isometrically on an $\R$--tree with virtually cyclic edge
    stabilizers and without global fixed point.
    \item[(2)] Let $G $ be a finitely generated hyperbolic group with Kazhdan property (T). Then $Out(G )$ is
    finite.
\end{itemize}
\end{theorem}

Statement (2) follows immediately from (1) because property (T)
implies that every action by isometries on a real tree has a
global fixed point.

\begin{remark}
Theorem \ref{out}, (1), together with \cite[Theorem 9.5]{BF} imply
that if $G$ is a hyperbolic group and if $Out(G )$ is infinite
then either $G$ splits as an amalgamated product or as an HNN
extension over a virtually cyclic subgroup, or $G$ is itself
virtually cyclic.
\end{remark}

\noindent \textit{Examples of hyperbolic groups with property
(T)}:

\begin{itemize}
    \item uniform lattices of isometries of
    $\hip^n_{\mathbb{H}}$, $n\geq3$;
    \item all their hyperbolic quotients. The quotient of a group $G$ with property (T) also has property (T). But hyperbolicity is
    not automatically inherited by a quotient. Nevertheless, it
    appears that ``almost every'' quotient of a hyperbolic group
    is hyperbolic, in the following sense. A quotient of the group $G$ means the prescription of new
     relations. A different way of saying it is that, given some finite generating set $S$ of the group $G$, one chooses
      a set of reduced words in the alphabet $S$ and puts the condition that they become equal to $1$.
       Since we want a hyperbolic quotient we prescribe finitely many new relations, that is we pick finitely many reduced words in
       $S$. There are several ways to introduce the probabilistic language
       into the picture. One of them is as follows. Choose randomly $e^{\beta \ell }$ new
     relations among the reduced words of length $\ell $ in $S$. Given a certain property (*), count the number
     $N_{\beta , \ell }$ of choices that give a quotient
      with property (*). The \textit{probability that the quotient has property (*)} is the
       limit as $\ell\to \infty$ of the ratio of $N_{\beta , \ell }$
      over the number of all possible choices of  relations under the parameters given above.
        According to the results in \cite{Gr4}, \cite{Gr5} and \cite{Oll}, for every non-elementary
hyperbolic group $G$ and finite generating set $S$ of it,
     there exists $\alpha =\alpha (G ,S)>0$
     such that the following holds:
    \begin{itemize}
        \item for every $\beta < \alpha $, the probability that the quotient is non-elementary hyperbolic is $1$;
        \item if $\beta > \alpha $ then with probability $1$
        the quotient is either trivial or $\Z / 2\Z$.
    \end{itemize}

    In the particular case when $G$ is the free group of
    rank  $m$, $\free_m$, and $S$ is the set of $2m$
     generators, $\alpha =\frac{\ln (2m-1)}{2}$.

 This gives a large choice of relations and of hyperbolic quotients for any hyperbolic group $G$.
 It allows in particular for the possibility of constructing approximate copies
  of very complicated graphs in the Cayley graph of a quotient of $G$,
  without loosing the property of hyperbolicity.
     \item In \cite{Z}, a slightly different notion of random group is considered. Given $\free_m$ the free group
     of rank $m$, one chooses randomly relations of
     length $3$ - on the whole there are $2m(2m-1)^2$ candidates. Suppose that one chooses randomly
     $(2m-1)^{3\beta }$ relations. The \textit{probability that the quotient has
property (*)} is in this case the
       limit as $m\to \infty$ of the ratio of the number of choices of $(2m-1)^{3\beta }$ relations
       that give a quotient of $\free_m$ with property (*),
      over the number of all possible choices of relations.

     In\cite{Z} it is proved that:
     \begin{itemize}
        \item if $\beta >\frac{1}{3}$ then the probability
     that the quotient has property (T)
     is 1;
        \item if $\beta <\frac{1}{2}$ then with probability $1$ the quotient is non-elementary hyperbolic;
       \item consequently for $\beta \in \left(\frac{1}{3},\frac{1}{2}
       \right)$, with probability $1$ the quotient is both
       hyperbolic and with property (T).
     \end{itemize}
\end{itemize}

\bi

In the case of relatively hyperbolic groups, the following two
theorems provide a sample of results on the group of outer
automorphisms.

\begin{theorem}[\cite{DS}]\label{outrh}
Let $G$ be a group relatively hyperbolic with respect to $\{
H_1,\dots ,H_n \}$ and suppose that all $H_i$ are asymptotically
without cut-points. Then for every $i\in \{ 1,\dots ,n \}$, there
exists a homomorphism from a subgroup of index at most $n!$ in
$Out(G)$ to $Out(H_i)$.
\end{theorem}

A result more in the line of Theorem \ref{out} is the following.

\begin{theorem}[Theorem 1.7 in \cite{DStree}]\label{outrh2}
Let $G$ be a group relatively hyperbolic with respect to $\{
H_1,\dots ,H_n \}$ and suppose that no $H_i$ is relatively
hyperbolic with respect to proper subgroups. If $Out(G)$ is
infinite then one of the followings cases occurs:
\begin{itemize}
    \item[(a)] $G$ splits as an amalgamated product or
    HNN extension over a virtually cyclic subgroup;
    \item[(b)] $G$ splits as an amalgamated product or
    HNN extension over a parabolic subgroup.
\end{itemize}
\end{theorem}

\begin{remark}
In particular a group $G$ as in Theorem \ref{outrh2} which
moreover has property (T) has finite $Out(G)$.
\end{remark}

 \me

 \section{Groups asymptotically with(out) cut-points}

First we return to the list of examples of groups asymptotically
without
 cut-points, drawn after Remark \ref{roneend}, and discuss Examples 3 and 4.

\subsection{Groups with elements of infinite order in the center, not virtually cyclic}


Let us see what happens if the asymptotic cone $\co{G;1,d}$ of an
infinite
 group $G$
 has a cut-point. Proposition \ref{cutp} implies that it is a
 tree-graded space with respect to a set of pieces $\pp$ such that each piece is either a point or a geodesic subset
 without cut-point. In particular, if all pieces are points the cone is a
 tree. Note that by homogeneity in this case it can be either a
 line or a tree in which every point is a branching point.

 The case when one asymptotic cone is a line turns out to be quite
 particular.
\begin{proposition}[Corollary 6.2 in \cite{DS}]
A finitely generated group such that one asymptotic cone is a
point or a line is virtually cyclic.
 \end{proposition}

Let now $G$ be a non-virtually cyclic group with a central
infinite cyclic subgroup $\langle h \rangle$. We have to show that
$G$ cannot have cut-points in any
 asymptotic cone. Suppose that one of its asymptotic cones  $\co{G;1,d}$
 has cut-points. It follows that it is tree-graded and that it is
 not a line.

Every element $\zeta $ in the center of $G$ is an isometry with
the property that
  every $g\in G$ is translated by $\zeta$ at a fixed distance, as $\dist (\zeta g, g)=\dist (\zeta ,1)$. It
  is not difficult to deduce from this the following. For every $\epsilon
  >0$ there exists an isometry $h_\omega $ of $\co{G;1,d}$, $h_\omega \in G^\omega
  $, such that for every $x\in \co{G;1,d}$, $\dist (h_\omega (x)\, ,\,
  x)=\epsilon$.

  On the other hand no tree-graded space different
  from a line admits such set of isometries. This is clearly seen
  for instance in the particular case of a real tree with at least
  one branching point. There is no way in which to translate a small tripod
  in this tree such that all its points move at the same
  distance.

\subsection{Groups satisfying an identity, not virtually
cyclic}\label{id}

Again we argue by contradiction and suppose that such a group $G$
has an asymptotic cone $\co{G;1,d}$ with cut-points, consequently
an asymptotic cone which is tree-graded and different from a line.
By the argument in the end of Section \ref{prel}, the group
$G^\omega$ acts transitively on $\co{G;1,d}$. Note that the
$\omega$--ultrapower $\Pi_\omega G$ and its subgroup $G^\omega$
satisfy the same identity as $G$. Even more can be said:

\begin{lemma}[Lemma 6.15, \cite{DS}]\label{fsg}
Let $\omega$ be any ultrafilter. The group $G$ satisfies a law if
and only if its $\omega$--ultrapower $\Pi_\omega G$ does not
contain free non-Abelian subgroups.
\end{lemma}

If $\co{G;1,d}$ is a tree then $G^\omega$ cannot act on it by
fixing a point in the boundary of the tree \cite[$\S 6$]{DS}. This
fact and \cite[Proposition 3.7, page 111]{Chis} imply that
$G^\omega$ contains a free non-Abelian subgroup. This contradicts
Lemma \ref{fsg}.

So in what follows we may assume that $\co{G;1,d}$ is not a tree,
consequently that it contains at least one (hence by homogeneity
continuously many) pieces without cut-points which are not
singletons.

To conclude we need the following result.

\begin{proposition}[\cite{DS}]\label{gomega}
 Let $\free $ be a tree-graded space with at least one non-singleton
 piece, and let $\ggg$ be a
group acting transitively on $\free$ and permuting pieces. The
group $\ggg$ contains a non-Abelian free subgroup.
\end{proposition}

\noindent \textit{Outline of proof.} Let $\pp$ be the set of
pieces of $\free$, containing at least one piece $P$ different
from a point. It follows that $P$ has cardinality  $2^{\aleph_0}$.
Then it can be shown that for every pair of distinct points $a,b$
in $P$ there exists an isometry $g\in \ggg$ such that $g(P)\cap
P=\emptyset$, $g(P)$ projects onto $P$ in $a$ and $P$ projects
onto $g(P)$ in $g(b)$. We denote by $\Pi_x$ the set of all points
projecting in $P$ in the point $x$. Property 6 of tree-graded
spaces implies that for every $x\neq b$, $g\left( \Pi_x
\right)\subset \Pi_a$.

From this one can easily deduce that $g\iv (P)\cap P=\emptyset$,
$g\iv(P)$ projects onto $P$ in $b$ and $P$ projects onto $g\iv(P)$
in $g\iv(a)$. Moreover, for every $y\neq a$, $g\iv \left( \Pi_y
\right)\subset \Pi_b$.

Now we choose a second pair of distinct points $c,d$ in
$P\setminus \{ a,b \}$. We choose for this pair a second isometry
$h$ with the same properties as $g$ for $a,b$.

\begin{figure}[!ht]
\centering
\includegraphics{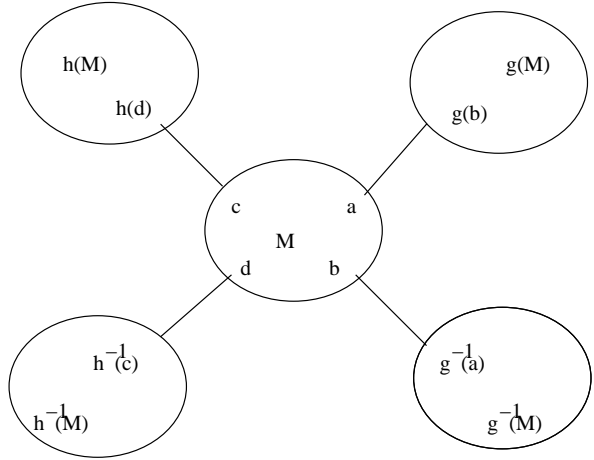}
\caption{The ping-pong argument.}
\label{fig8}
\end{figure}

It is not difficult to show by a ping-pong argument that $g$ and
$h$ generate a free group.\hspace*{\fill} $\Box$

\bigskip

As already mentioned, it turns out that uniformly amenable groups
are a particular case of groups satisfying a law. This observation
is maybe worth some explanations. Recall the following results.

\begin{theorem}[\cite{Wys}]\label{ultra}
Let $G$ be a countable discrete group.
\begin{itemize}
    \item[(1)] If $G$ is uniformly amenable then for any ultrafilter $\omega $ the ultrapower $\Pi_\omega G$
 is uniformly amenable.
    \item[(2)] If there exists an ultrafilter $\omega $ such that the ultrapower $\Pi_\omega G$ is
amenable then $G$ is uniformly amenable.
\end{itemize}
\end{theorem}

In particular, if the ultrapower of a discrete countable group is
amenable then it is uniformly amenable.

Now we recall some classical results.

\begin{proposition}
A subgroup $S$ of an amenable group $G$ is amenable.
\end{proposition}

\begin{remark}
Note that no other assumption is made on $S$ or $G$ - except the
amenability for $G$, as defined page 15 - not discreteness, nor
local compactness nor anything.
\end{remark}

\proof Take $\epsilon
>0$ arbitrary small. Take $K$ a finite subset in $S$. There exists a subset $F$ in $G$ such
that $\mathrm{card}\, K F < (1+\epsilon )\mathrm{card}\, F$.
Consider a graph whose vertices are the elements of the set $F$,
and whose edges correspond to the pairs of points $(f_1,f_2)\in
F\times F$ such that $f_2=kf_1$, where $k\in K$. Let $C$ be a
connected component of this graph with set of vertices
$\mathcal{V}_C$. Then $K\mathcal{V}_C$ does not intersect the sets
of vertices of other connected components. Hence there exists a
connected component $C$ such that $\mathrm{card}\, K\mathcal{V}_C
<(1+\epsilon )\mathrm{card}\, \mathcal{V}_C$ (otherwise if all
these inequalities have to be reversed, the sum of them gives a
contradiction with the choice of $F$). Without loss of generality,
we can assume that $\mathcal{V}_C$ contains $1$. Otherwise we can
shift it to $1$ by multiplying on the right by $c^{-1}$ for some
$c \in \mathcal{V}_C$. Then $\mathcal{V}_C$ can be identified with
a finite subset of $S$. Therefore $S$ contains a subset
$\mathcal{V}_C$ such that $\mathrm{card}\, K\mathcal{V}_C <
(1+\epsilon )\mathrm{card}\, \mathcal{V}_C$.\endproof

\begin{proposition}
A non-Abelian free group is not amenable.
\end{proposition}

A nice proof of this can be found in \cite[$\S 6.C$]{GLP}.

\begin{cor}\label{cna}
A group having a non-Abelian free subgroup is not amenable.
\end{cor}

\begin{remark}\label{rvn}
The observation that the existence of a free subgroup excludes
amenability was first made by J. von Neumann in \cite{vN}, the
very paper in which he introduced the notion of amenable group,
under the name of measurable group. It is this observation that
raised the question known later as \textit{the von Neumann
problem}: whether any non-amenable group contains a free
non-Abelian subgroup. In \cite{Ti} it was shown that for linear
groups the von Neumann problem has an affirmative answer, moreover
a linear group without any free non-Abelian subgroup is
solvable-by-finite. The first examples of non-amenable groups with
no (non-Abelian) free subgroups were given in \cite{Olsh}. In
\cite{Ad2} it was shown that the free Burnside groups $B(n,m)$
with $n\geq 2$ and $m\geq 665$, $m$ odd, are also non-amenable.
The first finitely presented examples of non-amenable groups with
no (non-Abelian) free subgroups were given in \cite{OlS}.
\end{remark}

Theorem \ref{ultra}, (1), Corollary \ref{cna} and Lemma \ref{fsg}
imply the following.

\begin{cor}[Corollary 5.9 in \cite{Kel}, Corollary 6.16 in \cite{DS}]\label{amen}
A finitely generated group which is uniformly amenable satisfies a
law. In particular no asymptotic cone of it has a cut-point.
\end{cor}

\subsection{Existence of cut-points in asymptotic cones and relative
hyperbolicity}\label{scutnrh}

A natural question to ask is the following.

\begin{ques}\label{qchar}
Can one improve the characterization of relatively hyperbolic
groups by their asymptotic cones given in Theorem \ref{thrhc} to:
a group $G$ is relatively hyperbolic if and only if all its
asymptotic cones have cut-points ?
\end{ques}

The``only if'' part is already proved. Concerning the ``if'' part
let us note that if an arbitrary asymptotic cone $\mathcal{C}$ of
$G$ has a cut-point then Proposition \ref{cutp} implies that it is
tree-graded with respect to some collection of pieces
$\mathcal{P}$ (which are either points or without cut-point).
Still it is not granted that there exists a finite set of
subgroups of $G$ such that all pieces are limit sets of left
cosets of these subgroups.

It turns out that the answer to Question \ref{qchar} is negative,
and that the property of having cut-points in every asymptotic
cone appears oftener than relative hyperbolicity. Here are some
examples of groups that are not relatively hyperbolic and have
cut-points in every asymptotic cone:

\begin{enumerate}
    \item  The mapping class group of an orientable finite type surface $S$ with
 $$3\cdot \, \mathrm{genus}(S) + \#\, \mathrm{punctures} \geq 5;$$
 The fact that it has cut-points
in any asymptotic cone is proved in \cite{Behrstock}. The fact
that it is not relatively hyperbolic can be deduced from arguments
in \cite{Bowh3m} and \cite{KN}, and it is explicitly proved in
\cite{AAS} and \cite{BDM}.
    \item Many right angled Artin groups \cite{BDM}.
    \item Fundamental groups of graph manifolds. They are not
    relatively hyperbolic according to \cite{BDM}, while they have
    cut-points in any asymptotic cone by arguments in \cite{KaL3}
    and \cite{KKL}.
\end{enumerate}

There also exists a metric example of the same sort, in which the
relative hyperbolicity is to be taken in its purely metric sense
given in \cite{DS}. More precisely, for any surface $S$ with
$3\cdot \, \mathrm{genus}(S) + \#\, \mathrm{punctures} \geq 9$,
the Teichm\"{u}ller space with the Weil-Petersson metric is not
relatively hyperbolic \cite{BDM}, while it has cut-points in any
asymptotic cone \cite{Behrstock}.

It follows from arguments in \cite{KKL} that the property of
having cut-points in all asymptotic cones is common to many
fundamental groups of non-positively curved compact manifolds.

\begin{proposition}[\cite{KKL}]
If $M$ is a compact non-positively curved manifold then either the
universal cover of $M$ is a symmetric space or $\pi_1(M)$ has
cut-points in any asymptotic cone of it.
\end{proposition}

\section{Open questions}\label{soq}

\begin{ques}
How does weak relative hyperbolicity behave with respect to
quasi-isometries ?
\end{ques}

The methods used for (strong) relative hyperbolicity no longer
work. Theorem \ref{s} again does not hold as can be easily seen by
taking $G=\Z^n$, $H=\Z^{n-1}\times \{ 0\}$ and $S=\Z^{n-1}$. A
quasi-isometric embedding of $S$ has no reason to stay close to a
left coset of $H$, as illustrated by many examples: it can be
transversal to all left cosets of $H$ or it can be composed of
many horizontal and vertical pieces etc.

Up to now there is no general result on the behavior up to
quasi-isometry of weakly relatively hyperbolic groups. In
\cite{KaL2}, \cite{Papasoglu:Zsplittings}, \cite{Dunwoody:Sageev},
\cite{MSW:QTOne} and \cite{MSW:QTTwo} strong quasi-isometric
rigidity results are proved for some particular cases of weakly
relatively hyperbolic groups---in fact all of them fundamental
groups of some graphs of groups. The notion of thick group
introduced in \cite{BDM} can be seen as a first attempt towards a
study of weakly relatively hyperbolic groups from the
quasi-isometry rigidity viewpoint.

\begin{ques}[``accessibility'' for relatively hyperbolic
groups]\label{last} Under which conditions does the process
described  in Corollary \ref{cor} have a terminal point, that is:
when does a relatively hyperbolic group $G$ have a list of
peripheral subgroups that are not relatively hyperbolic ? Does
this hold when $G$ is torsion-free, when it is finitely presented
? Note that both conditions are not satisfied by the inaccessible
groups of Dunwoody (see Remark \ref{coment}).
\end{ques}

We recall the standard theory of accessibility of groups, to which
this question relates. By Stalling's Ends Theorem \cite{Sta}, a
finitely generated group with more than one end splits as a free
product or HNN-extension with finite amalgamation. The question is
whether in an arbitrary finitely generated group one can keep on
doing this splitting until no more splitting is possible, that is
until all the factor groups are finite or one-ended. The answer is
positive for finitely generated torsion-free groups (the
Grushko-Neumann theorem) and for finitely presented groups
\cite{Du1}. But it is not true for all finitely generated groups
\cite{Du2}.

\begin{ques}
Given a group $G$ hyperbolic relative to the subgroups $H_1,\dots
, H_m$ can one say that the group $G$ has all asymptotic cones
isometric to each other, under the obvious necessary condition
that each $H_i$ has all asymptotic cones isometric to each other ?
\end{ques}

This would generalize the result of \cite{DP} from hyperbolic
groups to relatively hyperbolic groups.

\begin{ques}
Do relatively hyperbolic groups have uniform exponential growth
?\footnote{A positive answer to this question has been given in
\cite{Xie}.}
\end{ques}

A finitely generated group is said to have \textit{exponential
growth} if for some set of generators $S$ (hence for every $S$),
the growth function $B_S(n)=\mathrm{card}\, B(1,n)$ is
exponential. One can define $\alpha_S=\lim_{n\to \infty }
\frac{\ln B_S(n)}{\ln n}$ and then exponential growth means that
$\alpha_S >0$. One can also define $\alpha =\inf_{S}\alpha_S$. If
$\alpha
>0$ then the group is said to have \textit{uniform exponential
growth}. Hyperbolic groups for instance have uniform exponential
growth \cite{Kou}. There are also examples of groups having
exponential growth but not uniform exponential growth \cite{Wi}.
For a survey of the subject see \cite{dH}.

 The usual way in which uniform exponential growth is proved is to
 show that there exists some $n_0$ such that $B_S(1,n_0)$ contains
 two elements generating a free subgroup (or a free sub-semigroup), for every generating set
 $S$.

\me

\begin{ques}\label{q0}
 Is it true that an amenable group has at least one asymptotic cone without
 cut-points ?
\end{ques}

A stronger version of  this question was formulated by B. Kleiner:
do amenable groups have all asymptotic cones without cut-points
(that is, can the conclusion of Corollary \ref{amen} be extended
from uniformly amenable groups to amenable groups) ? In \cite{OOS}
Kleiner's question is answered in the negative: an example of an
amenable (and even elementary amenable) group with one asymptotic
cone a tree is constructed. Still, Question \ref{q0} remains open.




\section{Dictionary}\label{dic}

\begin{itemize}
    \item \textbf{Boundary at infinity.} Given $X$ either a simply connected Riemannian manifold of non-positive
    curvature (or more generally a $CAT(0)$--space) or an infinite graph, its
    boundary at infinity $\partial_\infty X$ is the quotient $\mathcal{R}/\sim
    $ of the set $\mathcal{R}$ of geodesic rays in $X$ with
    respect to the equivalence relation $r_1\sim r_2\Leftrightarrow \dist_H (r_1,r_2)<+\infty
    $.
    \item \textbf{(Abstractly) commensurable groups.} Two discrete groups $G_1$ and $G_2$ are called
\textit{abstractly commensurable} if they have finite index
subgroups that are isomorphic.

\item \textbf{Commensurable groups in an ambient larger group.} When both $G_1$ and $G_2$ are subgroups in a group $G$, we say that
$G_1$ and $G_2$ are \textit{commensurable (in $G$)} if there
exists $g\in G$ such that $G_1^g \cap G_2$ has finite index both
in $G_1^g$ and in $ G_2$.
\item \textbf{Commensurator.} In the above case the set of $g\in G$ such that $G_1^g \cap G_2$ has finite index
both in $G_1^g$ and in $ G_2$ is called \textit{the commensurator
of $G_1$ to $G_2$ in $G$}, and it is denoted $Comm_G(G_1,G_2)$.
When $G_1=G_2$ we simply write $Comm_G(G_1)$. Also, when there is
no possibility of confusion, we drop the index $G$.
\item \textbf{Convergence group.} It is a subgroup $G$ of Homeo$(\sph^1
)$ such that every sequence of distinct elements in $G$ contains a
subsequence $(g_n)$ for which there exist $x,y\in \sph^1$ with the
property that on $\sph^1\setminus \{ x,y \}$, $g_n$ converges to
$x$ and $g_n\iv $ converges to $y$ uniformly on compact subsets.
\item \textbf{Filter.} A \textit{filter} $\mathcal{F}$ over a set $I$ is a collection of subsets of $I$ satisfying the
following conditions:

\begin{itemize}
\item[$(F_1)$] If $A\in \mathcal{F}$, $A\subseteq B\subseteq I$, then $B\in
\mathcal{F}$;
\item[$(F_2)$] If $A, B\in \mathcal{F}$ then $A\cap B\in \mathcal{F}$;
\item[$(F_3)$] $\emptyset \not \in \mathcal{F}$.
\end{itemize}

For instance, if $I=\N$, the collection of all complementaries of
finite sets is a filter over $\N$, called \textit{the Fr\'echet
filter}.
\item \textbf{Fuchsian group.} It is a discrete subgroup of $PSL (2,\R )=\mathrm{Isom }(\hip^2_\R
)$.
\item \textbf{Fully residually * group (also called $\omega$--residually * group).} Here * represents a family of groups (finite groups, free groups etc.) A group $G$ is
fully residually * if for every finite subset $F$ in $G$ there
exists a homomorphism from $G$ onto a * group which is injective
on $F$.
\item \textbf{(Global) cut-point.} A point $p$ in a topological
space $X$ such that $X\setminus \{ p\}$ has several connected
components.
\item \textbf{Geodesic metric space (see Length metric
space).} It is a length metric space such that for every pair of
points, the shortest path joining them exists. By Hopf-Rinow
Theorem \cite{GLP} a complete locally compact length metric space
is geodesic.
\item \textbf{Hausdorff distance.} If $A$ and $B$ are two subsets
in a metric space $X$, then the Hausdorff distance $\dist_H (A,B)$
between $A$ and $B$ is the minimum of all $\delta >0$ such that
$A$ is contained in the $\delta$-tubular neighborhood of $B$ and
$B$ is contained in the $\delta$-tubular neighborhood of $A$. If
no such finite $\delta $ exists, one puts $\dist_H (A,B)=+\infty
$.

    \item \textbf{Hawaiian earring}. It is the topological space
    $\bigcup_{n\in \N } C\left( \left(0, \frac{1}{n} \right), \frac{1}{n}
    \right)$ with the topology induced from $\R^2$, where $C\left( \left(0,\frac{1}{n} \right), \frac{1}{n}
    \right)$ denotes the circle of center $\left(0, \frac{1}{n}
    \right)$ and of radius $\frac{1}{n}$. Its fundamental group is
    uncountable and non-free \cite{DES}.
    \item \textbf{Horoball, horosphere}. Let $\varrho $  be
    a geodesic ray in a simply connected Riemannian manifold of non-positive
    curvature (more generally in a $CAT(0)$--space) $X$. It defines a
    point at infinity $\alpha \in \partial_\infty X$. The \textit{open horoball $Hbo(\varrho)$ determined
    by} $\varrho$ is the union of open balls $\bigcup_{t>0}B(\varrho (t),
    t)$. Its closure $Hb(\varrho)$ is the \textit{closed horoball determined
    by} $\varrho$, and its boundary $H(\varrho )$ is the \textit{horosphere determined
    by} $\varrho$.

    Note that if $\varrho_1\, ,\, \varrho_2$ are
    asymptotic rays then there exists $\kappa >0$ such that $\nn_\kappa
    (Hbo(\varrho_i))=Hbo(\varrho_j)$, where $\{i,j\}=\{ 1,2\}$.
    Thus, one horoball defines all the other horoballs determined by rays in
    the same asymptotic class. Therefore it makes sense to no
    longer specify the ray, but only the point at infinity $\alpha $ corresponding to it,
    and to speak about all horoballs corresponding to rays with
    the same point at infinity as \textit{horoballs of basepoint
    $\alpha$}. For details on this notion see \cite{BH}.

    \item \textbf{Length (or path) metric space.} A metric space $(X, \dist_\ell )$ such
    that for every pair $x,y$ in $X$, $ \dist_\ell (x,y)= $
      the infimum of the lengths of the paths joining $x$ and $y$.
      A priori the path realizing the infimum might not exist.

     Note that given a metric space $(X, \dist )$, one can define the length of curves in it.
     Consequently one can define a ``length metric'' $\dist_\ell $ on $X$. The problem is that in this case $\dist_\ell$
      might take the value $+\infty $, because in case $x$ and $y$ are not
joined by any path of finite length, or simply by any path, one
puts $\dist_\ell (x,y)=+\infty$.
    \item \textbf{Net.} A \textit{net} in a metric space $X$
    is a subset $N$ of $X$ which is\begin{itemize}
        \item $\delta$-\textit{separated} for some $\delta >0$: for every $n_1,n_2\in
        N$, $\dist (n_1,n_2)\geq \delta$;
        \item $\epsilon$\textit{-covering} for some $\epsilon >0$:
        $X\subset \nn_\epsilon (N)$.
    \end{itemize}

    When more precision is needed, $N$ is also called $(\delta , \epsilon
    )$--net.
    \item \textbf{Proper metric space}. A metric space with the property
    that all its closed balls are compact. Note that by the Hopf-Rinow
    Theorem \cite{GLP} every complete, locally compact length
    metric space is proper.

    \item \textbf{Rank of a symmetric space (of non-positive sectional curvature)}. The maximal $n\in
    \N$ such that the $n$-dimensional Euclidean space can be
    embedded isometrically as a totally geodesic submanifold in
    the symmetric space.

    \item \textbf{Rank one symmetric space}. A symmetric space of non-positive sectional curvature
     and of rank one; also called hyperbolic space. With
     one exception, rank one symmetric spaces can be described as follows. Given $\K = \R , \C
    $ or $\bH$, where $\bH$ is the field of the quaternions, consider $x\mapsto \bar
    x$ the standard involution on $\K$ (the identity on $\R$, the conjugation on $\C$
    and on $\bH$), and consider on $\K^{n+1}\times \K^{n+1}$ the bilinear form
     $$
L(x,y)= x_0 {\bar y}_0-x_1 {\bar y}_1-x_2 {\bar y}_2\cdots -x_n
{\bar y}_n\, .
     $$ Let $G$ be the connected component of the identity of the stabilizer of
     $L$ in $SL(n+1, \K)$. The quotient $G/K$ with $K$ a
     maximal compact subgroup in $G$ is the $n$-dimensional $\K$-hyperbolic space
     $\mathbb{H}^n_\K$. It can be identified with
     $$
\mathbb{D}_\K = \{ x\in \K^n \mid x_1{\bar x}_1+x_2{\bar
x}_2+\cdots +x_n{\bar x}_n<1\}
     $$ and it can be endowed with a Riemannian metric invariant with respect to the action of $G$ (see \cite[$\S 19$]{Mos} for
     details). Besides the above spaces, there exists one more hyperbolic
     space, the Cayley hyperbolic plane of which a complete
     description can be found in \cite[$\S 19$]{Mos}. The real hyperbolic spaces have constant negative
     sectional curvature, while all the other hyperbolic spaces have pinched negative sectional curvature.

    \item \textbf{Reduced words.} Given an alphabet $S= \{ a_1,\dots ,a_n,a_1\iv , \dots ,
    a_n\iv\}$, a word in the alphabet $S$ is called \textit{reduced} if it
    does not contain subwords of the form $a_i a_i\iv$ or $a_i\iv
    a_i$.

    \item \textbf{Symmetric space (see Rank of a symmetric space, Rank one symmetric space)}.
    A complete simply connected Riemannian
    manifold $X$ such that for every point $p$ the geodesic symmetry $\sigma_p$ fixing $p$ defined by
    $\sigma_p (\exp_p (v))= \exp_p(-v)$ for every $v\in T_p X$ is a global isometry of $X$.
    The connected component of the
    identity of the group of isometries of $X$, which we denote by $G$, acts transitively on
    $X$. Therefore $X$ is a homogeneous
    space and can be identified with a coset
    space $G/K$, where $K$ is the stabilizer of a point in $X$
    (and also a maximal compact subgroup of $G$). Details on the
    notion can be found in \cite{He}.

    \item \textbf{Tubular neighborhood}. For a set $A$ in a metric space $X$ and for $\delta >0$ we define
    the $\delta$-\textit{tubular neighborhood $\nn_\delta (A)$ of }$A$ as
    the set
$$
\{ x\mid \dist (x,A)< \delta \}\, .
$$
    \item \textbf{Ultrafilter (see Filter).} An \textit{ultrafilter} over a set
    $I$ is a filter $\mathcal{U}$ over $I$ which is a maximal element in the
    ordered set of all filters over $I$ with respect to the
    inclusion. An ultrafilter can also be defined as a collection of subsets of $I$ satisfying the
    conditions $(F_1),\, (F_2),\, (F_3)$ defining a filter and the
    additional condition:
\begin{center}
 $(F_4)$  For every $A\subseteq I$ either $A\in \mathcal{U}$ or $I\setminus
A\in \mathcal{U}$.
\end{center}

A \textit{non-principal ultrafilter} is an ultrafilter containing
the Fr\'echet filter.

    \item \textbf{Virtually *.} A group is said to have property *
    \textit{virtually} if a finite index subgroup has the property *.
\end{itemize}

\end{document}